\documentclass[11pt]{amsart}
\newtheorem{theoreme}{Theorem}[section]
\newtheorem{fait}[theoreme]{Fact}
\newtheorem{lemme}[theoreme]{Lemma}
\newtheorem{proposition}[theoreme]{Proposition}

\newtheorem{corollaire}[theoreme]{Corollary}

\input xy
\input{std_macros_new.sty}
\xyoption{all}
\usepackage[normalem]{ulem}
\usepackage{enumerate}
\usepackage[a4paper, left=2cm,right=1.5cm,bottom=2cm,top=2cm]{geometry}
\usepackage{color}
\usepackage{graphicx,epsfig}
\usepackage{amscd}
\usepackage{amssymb}
\usepackage{enumitem}
\usepackage{ulem}
\usepackage{hyperref}
\hypersetup{colorlinks,allcolors=blue}

  \newcommand{\Sdeg}{\hbox{\rm deg}}
  \newcommand{\Sdim}{\hbox{\rm dim}}

  \newcommand{\Sdet}{\hbox{\rm det}}
 \newcommand{\SH}{\hbox{\rm H}}
  \newcommand{\SGL}{\hbox{\rm GL}}

      \newcommand{\SHom}{\hbox{\rm Hom}}
           \newcommand{\SEnd}{\hbox{\rm End}}

\title{Geometric monodromy  --  semisimplicity and maximality}

\author{Anna Cadoret, Chun-Yin Hui and Akio Tamagawa}

\begin{document}
\maketitle
 
 \begin{abstract} Let $X$ be a connected scheme, smooth and separated  over an algebraically closed field $k$ of characteristic $p\geq 0$, let $f:Y\rightarrow X$ be a smooth proper morphism and $x$ a geometric point on $X$. We prove that the tensor invariants of bounded length $\leq d$ of $\pi_1(X,x)$ acting on the \'etale cohomology groups $\SH^*(Y_x,\F_\ell)$ are the reduction modulo-$\ell$ of those  of $\pi_1(X,x)$ acting on   $\SH^*(Y_x,\Z_\ell)$ for $\ell $ greater than a constant depending only on $f:Y\rightarrow X$, $d$. We apply this result to show that the geometric variant with $\F_\ell$-coefficients of the Grothendieck-Serre semisimplicity conjecture -- namely that $\pi_1(X,x)$ acts semisimply on $\SH^*(Y_x,\F_\ell)$ for $\ell\gg 0$ -- is equivalent to the condition that the  image of $\pi_1(X,x)$ acting on   $\SH^*(Y_x,\Q_\ell)$ is `almost maximal' (in a precise sense; what we call `almost hyperspecial') with respect to the group of $\Q_\ell$-points of its Zariski closure. Ultimately, we prove the geometric variant with $\F_\ell$-coefficients of the Grothendieck-Serre semisimplicity conjecture. \scriptsize{\begin{center} 2010 \textit{Mathematics Subject Classification.} Primary: 14F20, 20G25; Secondary: 20E28.\end{center}} 
\end{abstract}\textit{} \\

 \section{Introduction}\label{Sec:Introduction}
 Let $X$ be a connected scheme, smooth and separated  over an algebraically closed field $k$ of characteristic $p\geq 0$ and let $f:Y\rightarrow X$ be a smooth proper morphism. As all the objects are of finite type over the base,  we may assume that $ X=X_0\times_{k_0}k$ for some smooth, separated and geometrically connected scheme $X_0$ over a finitely generated subfield $k_0\subset k$   and that $f:Y\rightarrow X$ is the base-change over $k$ of a smooth proper morphism of $k_0$-schemes $f_0:Y_0\rightarrow X_0$. In the following, we  always use the notation $f_0:Y_0\rightarrow X_0 / k_0$ for such a model. By the smooth-proper base-change theorem, for every prime $\ell\not=p$, the higher-direct image sheaves $R^*f_{0*}\Z/\ell^n$ are locally constant constructible hence, for every  geometric point $x$ on $ X$, they give rise to  continuous actions of the \'etale fundamental group 
 $ \pi_1(X_0,x)$ on $(R^*f_{0*}\Z/\ell^n)_{x}\simeq \SH^*(Y_{x},\Z/\ell^n)$,  $(R^*f_{0*}\Z_\ell)_{x}\simeq \SH^*(Y_{x},\Z_\ell)$ and $(R^*f_{0*}\Q_\ell)_{x}\simeq \SH^*(Y_{x},\Q_\ell)$.\\
  
 \noindent The aim of this paper is to prove  the following  two  statements about the restriction of these representations to the geometric \'etale fundamental group $\pi_1(X,x)$ (note that they are independent of the model $f_0:Y_0\rightarrow X_0 / k_0$).

  \begin{theoreme}\label{Th:GSS}  The following holds.  
  \begin{itemize}[leftmargin=* ,parsep=0cm,itemsep=0cm,topsep=0cm]
 \item (\ref{Th:GSS}) For  $\ell\gg 0$ (depending on $f:Y\rightarrow X$), the action of $\pi_1(X,x)$ on $\SH^*(Y_x,\F_\ell)$ is semisimple.  \end{itemize}
 \end{theoreme}
 
\noindent The assertion (\ref{Th:GSS}) is the natural geometric variant with $\F_\ell$-coefficients  of the Grothendieck-Serre semisimplicity conjecture, stating that the action of $\pi_1(X_0,x)$  on $\SH^*(Y_x,\Q_\ell)$ is semisimple for $\ell\neq p$ (see for instance p. 109 of the version of Tate's Woods Hole talk in \cite{TateWH}; see also Section \ref{Conj:Tate}). The geometric variant with $\Q_\ell$-coefficients of this conjecture is a celebrated theorem of Deligne, proved in \cite{Weil2}  (see below for details).   \\
  \begin{theoreme}\label{Th:GSSMax} The assertion (\ref{Th:GSS}) is equivalent to the following.
  \begin{itemize}[leftmargin=* ,parsep=0cm,itemsep=0cm,topsep=0cm]
 \item (\ref{Th:GSSMax}) After replacing $X$ by a connected \'etale cover and for  $\ell\gg 0$ (depending on $f:Y\rightarrow X$), the image of $\pi_1(X,x)$ in the group of $\Q_\ell$-points of its Zariski closure  in $\SGL(\SH^*(Y_x,\Q_\ell))$ is almost hyperspecial.
 \end{itemize}

 \end{theoreme}

 \noindent  Given a connected semisimple group $\mathcal{G}$ over $\Q_\ell$, recall that a compact subgroup $\Pi\subset \mathcal{G}(\Q_\ell)$ is called hyperspecial if there exists a semisimple group scheme $\frak{G}$ over $\Z_\ell$ with generic fiber $\mathcal{G}$ and such that $\Pi=\frak{G}(\Z_\ell)$. When they exist, hyperspecial subgroups are the compact subgroups of maximal volume in $\mathcal{G}(\Q_\ell)$. We say that a compact subgroup $\Pi\subset \mathcal{G}(\Q_\ell)$ is almost hyperspecial if $(p^{sc})^{-1}(\Pi)\subset \mathcal{G}^{sc}(\Q_\ell)$ is hyperspecial, where $p^{sc}:  \mathcal{G}^{sc}\rightarrow  \mathcal{G}$ denotes the simply connected cover. We refer to the beginning of Section 8 for further details. \\
 
\noindent A motivation for the reformulation (\ref{Th:GSSMax}) of (\ref{Th:GSS}) in terms of maximality is its potential applications to problems requiring large monodromy assumptions (see the introduction of \cite{Hall}  and references therein for a survey of  applications of large monodromy results).\\

\noindent The assertion (\ref{Th:GSS}) was previously known   in the following cases
 \begin{itemize}[leftmargin=* ,parsep=0cm,itemsep=0cm,topsep=0cm]
\item (\ref{Th:GSS}.1) when $p=0$;
\item (\ref{Th:GSS}.2) if one replaces $\F_\ell$-coefficients with $\Q_\ell$-coefficients;  
\item (\ref{Th:GSS}.3) when $f:Y\rightarrow X$ is an abelian scheme (or more generally for $\SH^1(Y_x,\F_\ell)$) and $p$ is arbitrary \cite{Zarhin} or a family of K3-surfaces and $p\not= 2$ \cite{SZ}
\end{itemize}
 while  (\ref{Th:GSSMax}) was previously known when $p=0$ (see for instance \cite[Rem. 2.5]{Adelic}). Let us also point out that an arithmetic variant of (\ref{Th:GSSMax}) holds over a set of primes of density one and for arbitrary systems of compatible rational semisimple $\ell$-adic representations \cite{LarsenMax}.\\

\noindent Let us recall the proofs of (\ref{Th:GSS}.1), (\ref{Th:GSS}.2). For (\ref{Th:GSS}.1), we may assume $k\subset \C$ and $x\in X(\C)$. The topological fundamental group $\Pi:=\pi_1^{top}(X_{\C}^{an},x)$ acts on the singular cohomology $H:=\SH^*(Y_{x}^{an},\Z) $. As $H $ is a finitely generated $\Z$-module, $H\otimes\F_\ell\simeq \SH^* (Y_{x}^{an},\F_\ell)$ for $\ell\gg 0$. Thus, by comparison of singular/\'etale cohomology and of topological/\'etale fundamental group, the image of  $ \pi_1(X,x)$ acting on $  \SH^*(Y_{x},\F_\ell)$ identifies with the image of $\Pi$ acting on $H\otimes \F_\ell$  for $\ell\gg 0$. The fact that  $\Pi$ acts semisimply on $H\otimes \F_\ell$ for $\ell\gg 0$ then follows formally \cite[Lem. 2.5]{Genus1} from the Hodge-theoretical fact  that $\Pi$  acts semisimply on $H\otimes\Q$ \cite[Cor. 4.2.9 (a)]{DeligneHodge2}. For (\ref{Th:GSS}.2), by standard specialization arguments (see for instance Subsection \ref{Subsec:Specialization}) we may assume $k_0 $ is a finite field and $k=\overline{k}_0$. Set $\mathcal{F}:=R^wf_*\Q_\ell$ and consider the largest maximal semisimple  smooth subsheaf $\mathcal{F}'\subset \mathcal{F}$. As $\mathcal{F}'_x\subset \mathcal{F}_x$ is stable under the action of $\pi_1(X_0,x)$, the extension
 $$0\rightarrow \mathcal{F}'\rightarrow \mathcal{F}\rightarrow \mathcal{F}/\mathcal{F}'\rightarrow 0$$
 corresponds to a class in $\SH^1(X,\underline{\hbox{\rm Hom}}(\mathcal{F}/\mathcal{F}',\mathcal{F}'))^F,$
 where $F$ denotes the geometric Frobenius on $X$. As $\mathcal{F}$ is  smooth, pure of weight $w$, $\underline{\hbox{\rm Hom}}(\mathcal{F}/\mathcal{F}',\mathcal{F}')$ is  smooth, pure of weight $0$ hence $\SH^1(X,\underline{\hbox{\rm Hom}}(\mathcal{F}/\mathcal{F}',\mathcal{F}'))$ is mixed of weights $\geq 1$ and  $\SH^1(X,\underline{\hbox{\rm Hom}}(\mathcal{F}/\mathcal{F}',\mathcal{F}'))^F=0$ (see \cite[(3.4)]{Weil2}).
  \\
 
\noindent So, when $p=0$, the essential ingredient is comparison between complex and \'etale topology, which provides an underlying $\Z$-structure for the action of  $ \pi_1(X,x)$ on $\SH^*(Y_x,\Z_\ell) $ and 
reduces (\ref{Th:GSS}.1) to a Hodge-theoretical statement by reduction modulo-$\ell$ for $\ell\gg 0$. When $p>0$,  such  an underlying $\Z$-structure is no longer available. For $\Q_\ell$-coefficients,  we can resort to Deligne's theory of weights. But such a theory does not exist for $\F_\ell$-coefficients.  Our basic strategy  is to  combine both aspects, namely  try and deduce  (\ref{Th:GSS}) from (\ref{Th:GSS}.2) by a reduction modulo-$\ell$ argument which involves Deligne's weight theory, in particular, the following two consequences of it:\\

\begin{itemize}[leftmargin=* ,parsep=0cm,itemsep=0cm,topsep=0cm]
\item (Fact \ref{Gabber}) The $\SH^*(Y_{x},\Z_\ell)$ are torsion-free for $\ell\gg 0$;
\item (Fact \ref{Deligne}) The $\SH^*(Y_{x},\Q_\ell)$ ($\ell$: prime $\not=p$) form a compatible system of $\Q$-rational representations.\\
\end{itemize}
The  combination of Fact \ref{Gabber} and Fact \ref{Deligne} provides a weak replacement for the $\Z$-structure in characteristic $0$, allowing us to `glue together' the various representations with $\F_\ell$-coefficients by means of the reduction modulo-$\ell$ of the characteristic polynomials of Frobenii.\\
 
 \noindent These ideas have already been exploited to obtain structural results  about the image of  $ \pi_1(X,x)$ acting on $  \SH^*(Y_{x},\F_\ell)$. For instance: \\
 
 \begin{itemize}[leftmargin=* ,parsep=0cm,itemsep=0cm,topsep=0cm]
\item (Fact \ref{JNTI}) After possibly  replacing $X$ by a connected \'etale cover and for $\ell\gg 0$,  the image of $ \pi_1(X,x)$ acting on $  \SH^*(Y_{x},\F_\ell)$ is perfect and generated by its order-$\ell$ elements. \\
\end{itemize}
\noindent This seemingly technical statement also plays a crucial part in our arguments. \\

\noindent However, to achieve the proofs of Theorem \ref{Th:GSS}  and Theorem \ref{Th:GSSMax}, more information is required. The way to grasp the missing information is Tannakian: instead of considering only $\SH^*(Y_{x},\Z_\ell)$, we  consider all possible tensor constructions of bounded length built from $\SH^*(Y_{x},\Z_\ell)$. Behind this is the observation that the image of $ \pi_1(X,x)$ acting on $  \SH^*(Y_{x},\Q_\ell)$ is captured by its Zariski closure while the image of   $ \pi_1(X,x)$ acting on $  \SH^*(Y_{x},\F_\ell)$ is captured by its algebraic envelope in the sense of Nori (this is one place where Fact \ref{JNTI} is crucial). Both are algebraic groups hence should be reconstructible from their tensor invariants.  Whence the idea to  compare the tensor invariants for the action of  $ \pi_1(X,x)$ on   $  \SH^*(Y_{x},\Q_\ell)$ and  $  \SH^*(Y_{x},\F_\ell)$. This is the core result of our paper.  To state it, we  need some notation. \\

\noindent For every partition $\lambda$ of an integer $d\geq 0$ let $c_\lambda\in \Z[\mathcal{S}_d]$ denote the associated Young symmetrizer and write $n_\lambda:=\frac{d!}{d_\lambda}$, where $d_\lambda$ is the dimension of the irreducible representation of the symmetric group $\mathcal{S}_d$ defined by $c_\lambda$. Then $\frac{1}{n_\lambda}c_\lambda\in \Z[\frac{1}{n_\lambda}][\mathcal{S}_d]$ is an idempotent. Fix a prime $\ell$ such that $d<\ell$. Let $\Lambda_\ell$ denote $\Z_\ell$ or $\F_\ell$ and let $\Pi$ be a profinite group acting continuously on  a finitely generated free $\Lambda_\ell$-module $M$. Let $\mathcal{S}_d$ act on $M^{\otimes d}$ on the right; this action commutes with the one of $\Pi$ thus 
$$S_\lambda(M):=\frac{1}{n_\lambda}c_\lambda(M^{\otimes d})\subset M^{\otimes d}$$
again gives a representation of $\Pi$. For $\lambda,\lambda^\vee$ partitions of integers $d,d^\vee$ respectively, write 
$$S_{\lambda,\lambda^\vee}(M):=S_\lambda(M)\otimes S_{\lambda^\vee}(M^\vee),$$
where $(-)^\vee$ denotes the $\Lambda_\ell$-dual.\\

\noindent For a profinite group $\Pi$ (in practice, $\Pi$ will be $\pi_1(X,x)$) acting continuously on  a finitely generated free $\Z_\ell$-module $M$, consider the following equivalent properties
$$
(Inv,M)\;\;\begin{tabular}[t]{ll}
 (i)& $M^\Pi\otimes\F_\ell\simeq (M\otimes\F_\ell)^\Pi$;\\
 (ii)& $\SH^1(\Pi,M)[\ell]=0$;\\
(iii)&  $\SH^1(\Pi,M)$ is torsion-free.\\
\end{tabular}$$

\noindent Eventually, for a group $\Pi_0$ acting on a module $M$ and a morphism $\Pi\rightarrow \Pi_0$, let $M|_{\Pi}$ denote the $\Pi$-module obtained from $M$ by restriction of the action from $\Pi_0$ to $\Pi$.

  \begin{theoreme}\label{Th:TensorInvariants} For all integers $d,d^\vee\geq  0$, partitions $\lambda, \lambda^\vee$ of $d, d^\vee$ respectively  and $\ell\gg 0$ (depending on $f$, $d, d^\vee$)
  \begin{itemize}[leftmargin=* ,parsep=0cm,itemsep=0cm,topsep=0cm]
  \item (\ref{Th:TensorInvariants}.0)  the property $(Inv, S_{\lambda,\lambda^\vee}(\SH^*(Y_{x},\Z_\ell)))$ holds. 
  \end{itemize}
\noindent Furthermore  $(Inv,M|_{\pi_1(X,x)})$ holds for every $\pi_1(X_0,x)$-module $M$ which is  of one of the following forms
\begin{itemize}[leftmargin=* ,parsep=0cm,itemsep=0cm,topsep=0cm]
\item (\ref{Th:TensorInvariants}.1)  a torsion-free quotient of $S_{\lambda,\lambda^\vee}(\SH^*(Y_{x},\Z_\ell))$; 
\item (\ref{Th:TensorInvariants}.2)   a submodule of $  S_{\lambda,\lambda^\vee}(\SH^*(Y_{x},\Z_\ell))$ with torsion-free cokernel.
\end{itemize} 
  \end{theoreme}

\noindent Theorem \ref{Th:TensorInvariants} applies in particular to $\lambda=(d)$ \textit{i.e.} $c_\lambda=\sum_{\sigma\in \mathcal{S}_d}\sigma$, for which $S_\lambda(M)=S^d(M)$ and to $\lambda=(1,\cdots, 1)$ \textit{i.e.} $c_\lambda=\sum_{\sigma\in \mathcal{S}_d}\hbox{\rm sign}(\sigma)\sigma$, for which $S_\lambda(M)=\Lambda^d(M)$.\\
 
 \noindent To prove  Theorem \ref{Th:TensorInvariants}, we reduce to the case where $X_0$ is a curve and $k_0$ is finite (this uses Bertini's theorem, de Jong's alterations and specialization of tame \'etale fundamental group). This allows us to use  Deligne's theory of weights and  the  machinery of \'etale cohomology.\\

\noindent To obtain Theorem \ref{Th:GSSMax}, we follow the above rough Tannakian strategy. Theorem \ref{Th:TensorInvariants} enables us to show that, for $\ell\gg 0$, the Nori envelope of the image of $\pi_1(X,x)$ acting on $\SH^*(Y_x,\F_\ell)$ identifies with the reduction modulo-$\ell$ of the Zariski closure $\frak{G}_{\ell^\infty}$ of the image of $\pi_1(X,x)$ acting on $\SH^*(Y_x,\Z_\ell)$ (Theorem \ref{Cor:SSConsequence1}). This  in turn enables us to show that there exists a constant $C\geq 1$ (depending only on $f$) such that the image  of $\pi_1(X,x)$ acting on $\SH^*(Y_x,\Z_\ell)$ has index $\leq C$ in   $\frak{G}_{\ell^\infty}(\Z_\ell)$   (what we call weak maximality -- see (\ref{Cor:SSConsequence1}.2)). Using this, we can give several equivalent formulations of (\ref{Th:GSS}), among which are that   $\frak{G}_{\ell^\infty}$ is semisimple  (\ref{Cor:SSConsequence2}.4) and  (\ref{Th:GSSMax}) (see Corollary \ref{Cor:SSConsequence3}).\\
 
 \noindent The reformulation (\ref{Cor:SSConsequence2}.4) raises a general question: given a connected semisimple group $G$ over $\Q_\ell$ together with a faithful finite-dimensional $\Q_\ell$-representation $V$ and a lattice $H\subset V$, can one exploit tensor invariants data (as in Theorem \ref{Th:TensorInvariants}) to deduce that the Zariski closure $\frak{G}$ of $G$ in $\SGL_H$ is a semisimple model of $G$ over $\Z_\ell$? This question led to a first complete proof of Theorem \ref{Th:GSS} (Section \ref{Proof1}) -- which is entirely due to the second author, Chun-Yin Hui. More precisely, the key-result is that, under mild assumptions,  the semisimplicity of $\frak{G}$ is encoded by a finite explicit list of tensor invariants dimensions (Theorem \ref{Th:Lie}). This criterion is then applied to $H:=\SH^*(Y_{x},\Z_\ell)$ using Theorem \ref{Th:TensorInvariants} and the tools developed for the proof of Theorem \ref{Th:GSSMax}. Theorem \ref{Th:Lie} relies on Lie theory and is of independent interest.\\
 
 \noindent After this first proof of Theorem \ref{Th:GSS} was obtained, we completed a second proof (Section \ref{Proof2}) which is cohomological and reminiscent of  the argument of Deligne.   This second proof  requires Theorem \ref{Th:TensorInvariants}, Fact \ref{Gabber}, Fact \ref{Deligne} and Fact \ref{JNTI} but involves no additional group-theoretical machinery. \\  
 
 \noindent We conclude by observing that Theorem \ref{Th:TensorInvariants} and Theorem \ref{Th:GSS} imply that the positive characteristic variant of the (arithmetic) Grothendieck-Serre-Tate conjectures with $\F_\ell$-coefficients follow from the usual Grothendieck-Serre-Tate conjectures (arithmetic, with $\Q_\ell$-coefficients) (Corollary \ref{Final}).\\

\noindent The paper is divided into three parts. In Part I (Sections 1-6), we review the properties of \'etale cohomology involved in the proofs of our main results and establish Theorem \ref{Th:TensorInvariants}; here, Deligne's weight theory is ubiquitous. In Part II (Sections 7-8), we develop the group-theoretical machinery leading to the proof of Theorem \ref{Th:GSSMax}. Part III (Sections 9-11) is devoted to the proofs of Theorem \ref{Th:GSS} and to the application to the (arithmetic) Grothendieck-Serre-Tate conjectures with $\F_\ell$-coefficients.  The second proof of Theorem \ref{Th:GSS} (Section \ref{Proof2}) can be read just after Part I. \\

\noindent\textbf{Acknowledgments:} Anna Cadoret was partly supported by the ANR grant ANR-15-CE40-0002-01. Chun-Yin Hui was supported by the National Research Fund, Luxembourg and cofounder under the Marie Curie Actions of the European Commission (FP7-COFUND).  Akio Tamagawa was partly supported by JSPS KAKENHI Grant Numbers 22340006, 15H03609.  The authors  thank the unknown referees for their comments, which helped improve the exposition of the paper. They   also thank Brian Conrad, Philippe Gille, Wilberd van der Kallen, Tamas Szamuely and Jilong Tong for their interest and comments. \\

\begin{center} \textbf{\fontfamily{lmr}\selectfont{PART I: \'ETALE COHOMOLOGY}}
\end{center}

 \section{Notation, conventions}\label{Sec:Notation}
 \subsection{} Given a field $k_0$ and an algebraically closed field $k$ containing $k_0$, write $\pi_1(k_0,k)$ (or simply $\pi_1(k_0)$) for $\hbox{\rm Aut}(\overline{k}_0/k_0)$, where $\overline{k}_0$ denotes the separable closure of $k_0$ in $k$. From a scheme-theoretic point of view, writing   $x$ for the geometric point $\hbox{\rm Spec}(k)\rightarrow \hbox{\rm Spec}(k_0)$, one has $\pi_1(k_0,k)=\pi_1(\hbox{\rm Spec}(k_0),x)$, which justifies the notation.
 If $k_0$ is finite, let   $F_{k_0}\in  \pi_1(k_0,k)$ (or simply $F$ if $k_0$ is clear from the context) denote the geometric Frobenius.\\
 
\subsection{}\label{Subsec:Not2} Given a prime $\ell$ and a profinite group $\Pi$, let $\hbox{\rm Rep}_{\Z_\ell}(\Pi)$ denote the category of finitely generated $\Z_\ell$-modules endowed with a continuous action of $\Pi$. Let $X_0$ be a connected scheme and $x$ a geometric point on $X_0$. Then the fiber functor  $\mathcal{F} \rightarrow  \mathcal{F}_{x}$ induces an equivalence\footnote{As $X$ is connected, there always exists an \'etale path   between two geometric points on $X$. So this equivalence of categories is independent of $x$ in a canonical way up to fixing an \'etale path from $x$ to any other geometric point.} from the category   $\mathcal{S}(X_0,\Z_\ell)$ of smooth $\Z_\ell$-sheaves on $X_0$ to $\hbox{\rm Rep}_{\Z_\ell}(\pi_1(X_0,x))$. Thus, if $P$ is a property of objects in $\hbox{\rm Rep}_{\Z_\ell}(\pi_1(X_0,x))$ (\textit{e.g.},  torsion-free, torsion, irreducible, semisimple \textit{etc.}) we will say that  $\mathcal{F}\in \mathcal{S}(X_0,\Z_\ell)$  has $P$ if $\mathcal{F}_{x}$ has $P$. The same considerations apply to the corresponding $\Q_\ell$-categories. In the following, we will often implicitly identify smooth $\Z_\ell$- (resp. $\Q_\ell$-) sheaves on $X_0$ and finitely generated $\Z_\ell$- (resp. $\Q_\ell$-) modules endowed with a continuous action of $\pi_1(X_0,x)$.\\

 \noindent For every $x_0\in X_0$ and geometric point $x$ over $x_0$, consider the natural action of $\pi_1(x_0,x)$ on $\mathcal{F}_{x}$ defined as the composition 
 $$\rho_{x_0}:\pi_1(x_0,x)\stackrel{x_0}{\rightarrow}\pi_1(X_0,x) \rightarrow\SGL(\mathcal{F}_{x}):=\hbox{\rm Aut}_{\Z_\ell}(\mathcal{F}_{x}).$$
   
 \noindent  Assume $X_0$ is geometrically connected and of finite type over a field $k_0$.  Let $x_0\in X_0$ and fix a geometric point $x:\hbox{\rm Spec}(k)\rightarrow X$ over $x_0$; write   $X:=X_0\times_{k_0}k$. Then the sequence 
   $$  \pi_1(X,x)\rightarrow  \pi_1(X_0,x)\rightarrow \pi_1(k_0,k)\rightarrow 1$$
   is exact. The functor $\mathcal{F}\rightarrow \mathcal{F}_{x}^{\pi_1(X,x)}$ from $\mathcal{S}(X_0,\Z_\ell)$ to $\hbox{\rm Rep}_{\Z_\ell}(\pi_1(k_0,k))$
coincides with the global section functor $\SH^0(X,-):\mathcal{S}(X_0,\Z_\ell) \rightarrow\hbox{\rm Rep}_{\Z_\ell}(\pi_1(k_0,k))$ and 
   the action of $\pi_1(k_0,k)$ on $\SH^0(X,\mathcal{F})$ by `transport of structure'  identifies with the labelled arrow $(*)$ in the commutative diagram:
  $$\xymatrix{&\pi_1(k_0,k)\ar[dr]^{(*)}&\\
  \pi_1(x_0,x)\ar[r]\ar[ur]\ar@/_1pc/[rr]_{\rho_{x_0}} &\hbox{\rm coker}(\pi_1(X ,x)\rightarrow \pi_1(X_0,x))\ar[u]^\simeq\ar[r]&\SGL(\mathcal{F}_x^{\pi_1(X ,x)}),}$$
which also shows that  the restriction to $ \mathcal{F}_x^{\pi_1(X ,x)}$ of the action of $  \pi_1(x_0,x)$ on $\mathcal{F}_x$ factors through the canonical morphism $  \pi_1(x_0,x)\rightarrow \pi_1(k_0,k)$. In particular, if $x_0\in X_0(k_0)$, the induced action of $  \pi_1(x_0,x)$ on $\mathcal{F}_x^{\pi_1(X ,x)}$ is independent of $x_0$. \\

\noindent For $\sigma\in \pi_1(x_0,x)$, write 
 $$P_{\sigma,x_0}^{\mathcal{F}}:= \Sdet(TId-\rho_{x_0}(\sigma),\mathcal{F}_{x}\otimes\Q_\ell)\in \Q_\ell[T].$$
   If the residue field $k(x_0)$ at $x_0$ is a finite field, we will simply write  $F_{x_0}:= F_{k(x_0)}$ and $P^{\mathcal{F}}_{x_0}:=P^{\mathcal{F}}_{F_{x_0},x_0}$ (or even simply $P_{x_0}$ if $\mathcal{F}$ is clear from the context).\\

 \noindent Let $q$ be a power of a prime number  and $w\in \Z$. A $q$-Weil number of weight $w$ is an algebraic number $\alpha$ such that $|\iota(\alpha)|=q^{\frac{w}{2}}$ for every complex embedding $\iota:\overline{\Q}\hookrightarrow\C$. If $X_0$ is of finite type over $\Z$, following \cite[(1.2)]{Weil2}, a smooth $\Z_\ell$-sheaf $\mathcal{F}$ on $X_0$ is said to be  pure of weight $w$ (resp.  mixed) if for every closed points $x_0\in X_0[\frac{1}{\ell}]$ the roots of $P^{\mathcal{F}}_{x_0} $ are $|k(x_0)|$-Weil numbers of weight $w$ (resp. if $\mathcal{F}$ admits a filtration whose  successive quotients are  pure; the weights of the non-zero quotients are then called the weights of $\mathcal{F}$). A smooth $\Z_\ell$-sheaf $\mathcal{F}$ on $X_0$ is said to be $\Q$-rational if for every closed point $x_0\in X_0[\frac{1}{\ell}]$, $P^{\mathcal{F}}_{x_0}$ is in $ \Q[T]$. Given a  set $L$ of primes, a system $\mathcal{F}_{\ell}$, $\ell\in L$ of
  smooth $\Z_\ell$-sheaves on $X_0$ is said to be $\Q$-rational compatible if each of the $\mathcal{F}_{\ell}$ is $\Q$-rational and if for every closed point  $x_0\in X_0$ the polynomials $ P^{\mathcal{F}_{\ell}}_{x_0}\in  \Q[T]$ are independent of $\ell$ (for $\ell$ not equal to the residue characteristic of $x_0$).

\subsection{}Given a prime $\ell$ and $0\not=P=\sum_{n\geq 0}a_nT^n\in \Q[T]$, we define the reduction modulo-$\ell$  $\overline{P^\ell} $ of $P$ to be the reduction modulo-$\ell$ of $ a(P)P \in \Z[T]$, where 
  $$a(P)=\prod_{p}p^{-\hbox{\rm\tiny min}\lbrace v_p(a_n),\; n\geq 0\rbrace}.$$
  Here, the product is over all rational primes and $v_p:\Q\rightarrow \Z\cup\lbrace \infty\rbrace$ is the $p$-adic valuation. Given an $\F_\ell[T]$-module $M$ and $P\in\Q[T]$, we say that $M$ is killed by $P$ if it is killed by $\overline{P^\ell}$.

   \section{Some consequences of the Weil conjectures}\label{EC} From now on, we retain the notation and conventions of the introduction for $f:Y\rightarrow X/k$ and $f_0:Y_0\rightarrow X_0/k_0$.\\
   
   \noindent By the smooth-proper base-change theorem $R^*f_{0*}\Z_\ell$ is a smooth $\Z_\ell$-sheaf on $X_0$ and $(R^*f_{0*}\Z_\ell)_{x}=\SH^*(Y_x,\Z_\ell)$ for every geometric point $x$ on $X$. The following facts all rely on the Weil conjectures.

    \begin{fait}{\textnormal{(\hspace{-0.01cm}\cite{Gabber} (projective case), \cite[Rem. 3.1.5]{Orgogozo})}}\label{Gabber}  The  smooth $\Z_\ell$-sheaves $R^*f_*\Z_\ell$ are  torsion-free (of finite constant rank) for $\ell\gg 0$. In particular, 
$$\SH^*(Y_x,\Z_\ell)\otimes\F_\ell\tilde{\rightarrow} \SH^*(Y_x,\F_\ell),\;\;\textnormal{for}\;\ell\gg 0.$$
\end{fait}
    \begin{fait}{\textnormal{(\hspace{-0.01cm}\cite[Cor. 3.3.9]{Weil2})}}\label{Deligne} Assume $k_0$ is finite (so that $X_0$ is of finite type over $\Z$). Then $R^*f_{0*}\Z_\ell$ ($\ell$: prime $\not=p$) is a $\Q$-rational compatible system.\end{fait}
      \begin{fait}{\textnormal{(\hspace{-0.01cm}\cite[Cor. 3.4.3, Cor. 1.3.9]{Weil2}, \cite[Prop. 1.1]{LarsenPinkAV}, \cite[Prop. 2.2]{LarsenPinkAV})}}\label{DeligneLP}  After possibly replacing $X_0$ by a connected \'etale Galois cover,  the  Zariski closure of the image of $\pi_1(X,x)$ (resp. $\pi_1(X_0,x)$) acting on $\SH^*(Y_x,\Q_\ell)$ is connected semisimple (resp. connected) for all $\ell\not=p$.\end{fait}
  \begin{fait}{\textnormal{(\hspace{-0.01cm}\cite[Thm. 1.1]{JN}, \cite[Fact 5.1]{TI})}}\label{JNTI} 
After possibly replacing $X$ by a connected \'etale Galois cover and for $\ell\gg 0$ depending only on $f:Y\rightarrow X$,
the  image of $\pi_1(X,x)$ acting on $\SH^*(Y_x,\F_\ell)$ is perfect and generated by its order $\ell$ elements.\end{fait}

  \noindent  Fact \ref{Gabber} is used several times in the proof of Theorem \ref{Th:TensorInvariants}: to reduce the proof of Theorem \ref{Th:TensorInvariants} to the case where $Y\rightarrow X$ is the base change of a smooth proper morphism $Y_0\rightarrow X_0$ with $X_0$ a curve over a finite field $k_0$ (see Subsection \ref{Subsec:Localization}) and then, combined with Fact \ref{Deligne}, to compare the actions   of $\pi_1(X,x)$ on $\SH^*(Y_x,\Q_\ell)$ and $\SH^*(Y_x,\F_\ell)$ for $\ell\gg 0$ (see Subsection \ref{Sec:TensorInvariants}).  Fact \ref{DeligneLP} and Fact \ref{JNTI} are  used in the proof of  (\ref{Th:TensorInvariants}.2) (see Subsection \ref{ProofTI2}) and they also play a crucial part in the proofs of Theorem \ref{Th:GSSMax} (Part II) and Theorem \ref{Th:GSS} (Part III). 
   
\section{Preliminary reductions}\label{Sec:Reductions}   Consider the following assertions.
$$\begin{tabular}[t]{ll}
$(Inv,f_0)$&For  $\ell\gg 0$, $(Inv,M|_{\pi_1(X,x)})$ holds for   every $\pi_1(X_0,x)$-module $M$ which is a torsion-free quotient\\
&  of  $\SH^*(Y_{x},\Z_\ell)$ or 
  a submodule of $\SH^*(Y_{x},\Z_\ell)$ with torsion-free cokernel.\\
&\\
$(GSS,f)$& For $\ell\gg 0$, the action of $\pi_1(X,x)$ on $\SH^*(Y_x,\F_\ell)$ is semisimple.
\end{tabular}$$
Let $(Inv)$ (resp. $(GSS)$) denote $(Inv, f_0)$ (resp. $(GSS,f)$) for every  smooth proper morphism $f_0:Y_0\rightarrow X_0$. The aim of this section is to prove the following.
\begin{proposition}\label{Prop:Reductions}
Assume (Inv,$f_0$) (resp. (GSS,$f$)) holds when  $k_0$ is a finite field, $k=\overline{k}_0$ and $X_0$ is  a smooth, separated and geometrically connected curve over $k_0$. Then Theorem \ref{Th:TensorInvariants} (resp. (\ref{Th:GSS})) holds.
\end{proposition}

\subsection{Localization}\label{Subsec:Localization} 
 We begin with a formal group-theoretic lemma.
\begin{lemme}\label{Lem:Localization} Let   $\Pi$ be a profinite group and $U\subset \Pi$ a normal open subgroup. Fix a prime $\ell$ not dividing $[\Pi:U]$.  The following hold.
  \begin{itemize}[leftmargin=* ,parsep=0cm,itemsep=0cm,topsep=0cm]
 \item (\ref{Lem:Localization}.1)  Let $H$ be a finitely generated  free  $\Z_\ell$-module endowed with a continuous action of $\Pi$.  If  the  sequence $0\rightarrow H^U\stackrel{\ell}{\rightarrow} H^U\rightarrow (H\otimes\F_\ell)^U\rightarrow 0$  is exact then the sequence $0\rightarrow H^\Pi\stackrel{\ell}{\rightarrow} H^\Pi\rightarrow (H\otimes\F_\ell)^\Pi\rightarrow 0$  is also exact.
  \item (\ref{Lem:Localization}.2) Let $H$ be a finitely generated    $\F_\ell$-module endowed with a continuous action of $\Pi$. If  $U$ acts semisimply on $H$ then  $\Pi$ also acts semisimply on $H$.
  \end{itemize} 
\end{lemme}
\noindent\textit{Proof.} The assertion (\ref{Lem:Localization}.1) follows from the fact that $\SH^1(\Pi/U,H^U)[\ell]=0$ by applying the functor $(-)^{\Pi/U}$ to the short exact sequence $0\rightarrow H^U\stackrel{\ell}{\rightarrow} H^U\rightarrow (H\otimes\F_\ell)^U\rightarrow 0$. The assertion (\ref{Lem:Localization}.2) follows from the fact that for every $\Pi$-submodule $H'\subset H$, $U$-equivariant projector $p_U:H\twoheadrightarrow H'$ and system of representatives $\pi_1,\dots,\pi_r$ of $\Pi/U$ the map $p_\Pi:H\twoheadrightarrow H'$ defined by
 $$p_\Pi(h)=\frac{1}{r}\sum_{1\leq i\leq r}\pi_i p_U(\pi_i^{-1}h) $$
 is a   $\Pi$-equivariant projector. $\square$ \\

   \noindent  As a result, to prove Theorem  \ref{Th:TensorInvariants} (resp. (\ref{Th:GSS})) we may   base-change $f_0:Y_0\rightarrow X_0$ over $X'_0\rightarrow X_0$ for any morphism $X'_0\rightarrow X_0$ inducing a morphism $\pi_1(X'_0,x')\rightarrow \pi_1(X_0,x)$ with normal open image. This applies for instance to open immersions or connected \'etale Galois covers.
 \subsection{Specialization}\label{Subsec:Specialization} Let $\Lambda_\ell$ denote $\Z_\ell$ or $\F_\ell$.
 
 \begin{lemme}\label{Lem:Specialization} There exist a finite field $\tilde{k}_0 $, a smooth, separated and geometrically connected curve $\widetilde{X}_0$ over $\tilde{k}_0$, a smooth proper morphism $\widetilde{f}_0:\widetilde{Y}_0\rightarrow \widetilde{X}_0$, and a normal open subgroup $U \subset \pi_1(X,x)$, such that the following holds. Write $\tilde{k}=\overline{\tilde{k}}_0$ and $\widetilde{f}:\widetilde{Y}\rightarrow \widetilde{X}$ for the base-change of $\widetilde{f}_0:\widetilde{Y}_0\rightarrow \widetilde{X}_0$ over $\tilde{k}$. Then, for every geometric point $\tilde{x}$ on $\tilde{X}$ and $\ell\gg 0$, there is an isomorphism $\SH^*(\widetilde{Y}_{\widetilde{x}},\Lambda_\ell)\simeq \SH^*(Y_{ x},\Lambda_\ell)$   such that the image of $\pi_1(\widetilde{X},\widetilde{x})$ acting on $\SH^*(\widetilde{Y}_{\widetilde{x}},\Lambda_\ell)$ identifies with the image of $U$ acting on $\SH^*(Y_x,\Lambda_\ell)$. Furthermore we may assume that the action of
 $\pi_1(\widetilde{X}_0,\widetilde{x})$ on $\SH^*(\widetilde{Y}_{\widetilde{x}},\Lambda_\ell)$ factors through the tame \'etale fundamental group $\pi_1(\widetilde{X}_0,\widetilde{x})\rightarrow \pi_1^t(\widetilde{X}_0,\widetilde{x})$. 
\end{lemme}
\noindent\textit{Proof.} We proceed in two steps.\\
 \begin{itemize}[leftmargin=* ,parsep=0cm,itemsep=0cm,topsep=0cm]
 \item Reduction to dim$(X_0)=1$: We may assume $X_0$ has dimension $\geq 1$. Using  \cite[Thm. 6.10]{Jouanolou} (see \cite[Ex. 3.1]{JN} for details),
we can construct a finitely generated field extension $K_0$ of $k_0$ and a closed curve $C_0\subset X_0\times_{k_0}K_0 $ smooth, separated, geometrically connected over $K_0$, such that for every geometric point $c$ on $ C$ mapping to $x$ on $X$ the
induced morphism
$\pi_{1}(C ,c)\rightarrow\pi_{1}(X ,x)$ is surjective. (Note that, if $x$ is fixed, we cannot ensure that there exists a geometric point $c$ on $C$ mapping to $x$ but this is not a problem since, to prove Lemma \ref{Lem:Specialization}, we may replace $x$ by any other geometric point - see Footnote 1). Here, we write $K:=\overline{K_0k}$ and $C:=C_0\times_{K_0}K$. So the conclusion follows from the fact that the resulting representation 
  $\pi_{1}(C_0,c)\rightarrow\pi_{1}(X_0,x) \rightarrow\hbox{\rm GL}(\SH^*(Y_x,\Lambda_\ell))$
identifies  (modulo the isomorphism $\SH^*((Y\times_XC)_{c},\Lambda_\ell)\simeq \SH^*(Y _{x},\Lambda_\ell)$ given by the smooth-proper base-change theorem) with the representation
  $\pi_{1}(C_0,c)  \rightarrow\hbox{\rm GL}(\SH^*((Y\times_XC)_{c},\Lambda_\ell))$
  associated to
  the base-change $Y_0\times_{X_0}C_0\rightarrow C_0$.  \\

 \item Reduction to dim$(X_0)=1$, $|k_0|<+\infty$ and $k=\overline{k}_0$: From the above, we may assume that  $X_0$ is a smooth, separated, geometrically connected curve  over $k_0$. After enlarging $k_0$, we may assume it has a smooth compactification $X_0^{cpt}$ over $k_0$. From   
 de Jong's alteration theorem \cite[Prop. 6.3.2]{Berthelot}, for every $u\in  X_0^{cpt}\setminus X_0$ there exists an open subgroup 
$U_{u}$ of the inertia group $I_{u}\subset\pi_1(X_0,x)$ at $u$ such 
that the image of $U_{u}$ in 
$\hbox{\rm GL}(\SH^*(Y_x, \Q_\ell))$ 
is unipotent  for all $\ell\neq p$. 
The image of $U_{u}$ in 
$\hbox{\rm GL}(\SH^*(Y_x, \Z_\ell))$ is also a pro-$\ell$ group for $\ell\gg 0$ (Fact \ref{Gabber}). Hence, after replacing $X_0$ by a connected \'etale Galois cover (Subsection \ref{Subsec:Localization}), we may assume that $\pi_{1}(X_0,x) \rightarrow\hbox{\rm GL}(\SH^*(Y_x, \Z_\ell))$ factors through the tame \'etale fundamental group $\pi_{1}(X_0,x)\rightarrow \pi_{1}^t(X_0,x)$ and even  that $\pi_{1}(X_0,x) \rightarrow\hbox{\rm GL}(\SH^*(Y_x, \F_\ell))$ does (Fact \ref{Gabber}). Then by the standard specialization arguments (specialization of tame \'etale fundamental group, smooth-proper base-change for \'etale cohomology), we may assume that $k_0$ is finite and $k=\overline{k}_0$.  $\square$
 \end{itemize} 

\subsection{End of proof of Proposition \ref{Prop:Reductions}}\label{Subsec:End}   Lemma \ref{Lem:Localization} and Lemma \ref{Lem:Specialization} show that if  $(Inv,f_0)$ (resp. $(GSS,f)$) holds when $X_0$ is a smooth, separated, geometrically connected curve over a finite field $k_0$ and $k=\overline{k}_0$ then $(Inv)$ (resp. $(GSS)$) holds. We now explain why $(Inv)$ implies  Theorem \ref{Th:TensorInvariants}. \\
 
\noindent Consider the $d$-fold fiber product  $$f_0^{[d]}:Y_0^{[d]}:=Y_0\times_{X_0}\cdots\times_{X_0}Y_0\rightarrow X_0.$$
 By construction $(Y^{[d]})_x=  (Y_x)^{[d]}$ and, as $\SH^*(Y_x,\Z_\ell)$ is torsion-free for $\ell\gg 0$, the K\"{u}nneth formula (for both $\Z_\ell$- and $\F_\ell$-coefficients) shows that the horizontal arrows in the canonical commutative square of graded $\pi_1(X_0,x)$-modules
   $$\xymatrix{\SH^*((Y^{[d]})_x,\Z_\ell)\ar[r]^\simeq \ar[d]& \SH^*(Y_x,\Z_\ell)^{\otimes d}\ar[d]\\
    \SH^*((Y^{[d]})_x,\F_\ell)\ar[r]^\simeq& \SH^*(Y_x,\F_\ell)^{\otimes d}}$$  
are isomorphisms. This induces a commutative square, whose horizontal arrows are still isomorphisms
 $$\xymatrix{\SH^*((Y^{[d]})_x,\Z_\ell)^{\pi_1(X,x)}\otimes\F_\ell\ar[r]^\simeq \ar[d]& (\SH^*(Y_x,\Z_\ell)^{\otimes d})^{\pi_1(X,x)}\otimes\F_\ell\ar[d]\\
    \SH^*((Y^{[d]})_x,\F_\ell)^{\pi_1(X,x)}\ar[r]^\simeq& (\SH^*(Y_x,\F_\ell)^{\otimes d})^{\pi_1(X,x)}}$$
 On the other hand, as 
 $f_0^{[d]}:Y_0^{[d]}\rightarrow   X_0$ is a smooth proper morphism,  (Inv, $f_0^{[d]}$) implies that for $\ell\gg 0$ (depending on $f:Y\rightarrow X$, $d$) the left vertical arrow is an isomorphism as well.\\
   
 \noindent We now conclude the proof of  Theorem \ref{Th:TensorInvariants}. For $\ell\gg 0$ (depending on $f:Y\rightarrow X$, $d, d^\vee$), $\SH^*(Y_x,\Z_\ell)$ is torsion-free so Lemma \ref{Lem:Formal} below reduces the assertion of   Theorem \ref{Th:TensorInvariants} to the statement that  $(Inv, M|_{\pi_1(X,x)})$ holds for every $\pi_1(X_0,x)$-module $M$ which is a torsion-free quotient of $\SH^*(Y_x,\Z_\ell)^{\otimes d}\otimes  (\SH^*(Y_x,\Z_\ell)^{\vee})^{\otimes d^\vee}$ or a submodule of $  \SH^*(Y_x,\Z_\ell)^{\otimes d}\otimes  (\SH^*(Y_x,\Z_\ell)^{\vee})^{\otimes d^\vee}$ with torsion-free cokernel. But by Poincar\'e duality and the K\"{u}nneth formula,  
 $\SH^*(Y_x,\Z_\ell)^{\otimes d}\otimes  (\SH^*(Y_x,\Z_\ell)^{\vee})^{\otimes d^\vee}$
 is isomorphic to  $\SH^*(Y^{[d+d^\vee]}_x,\Z_\ell)$ (after suitable Tate twists) as $\pi_1(X_0,x)$-module. So the conclusion follows from $(Inv, f_0^{[d+d^\vee]})$. $\square$
 
 \begin{lemme}\label{Lem:Formal} Let $\Pi$ be a profinite group acting continuously on  a finitely generated torsion-free $\Z_\ell$-module $M$ and let $0\leq d, d^\vee<\ell$ be  integers. Then, for every pair of partitions $\lambda,\lambda^\vee$ of $d,d^\vee$ respectively,  $S_{\lambda,\lambda^\vee}(M)$ is a direct factor of  $M^{\otimes d}\otimes M^{\vee\otimes d^\vee}$  as a $\Pi$-module. 
 \end{lemme}
 
\noindent\textit{Proof.} Write $S_\lambda'(M):=(1-\frac{1}{n_\lambda}c_\lambda)(M^{\otimes d})\subset M^{\otimes d}$; this is again a $\Pi$-submodule and we have a direct sum decomposition 
$$M^{\otimes d}=S_\lambda(M)\oplus S_\lambda'(M).$$ 
Similarly we have $$M^{\vee\otimes d^\vee}=S_{\lambda^\vee}(M^\vee)\oplus S_{\lambda^\vee}'(M^\vee).$$
This implies that $ S_{\lambda,\lambda^\vee}(M) $  is a direct factor of $ M^{\otimes d}\otimes M^{\vee\otimes d^\vee}$. $\square$

\subsection{}\hspace{-0.2cm}Proposition \ref{Prop:Reductions} thus reduces the proof of Theorem \ref{Th:TensorInvariants} to the following special case.  Assume $X_0$ is a  smooth, separated and geometrically connected curve over a finite field $k_0$ and that $k=\overline{k}_0$. Then
 \begin{theoreme}\label{Th:TensorInvariantsHeart} For $\ell\gg 0$ (depending on $f:Y\rightarrow X$),
 \begin{itemize}[leftmargin=* ,parsep=0cm,itemsep=0cm,topsep=0cm]
 \item (\ref{Th:TensorInvariantsHeart}.0) (Inv, $\SH^*(Y_x,\Z_\ell)$) holds.
 \end{itemize}
 \noindent Furthermore  $(Inv,M|_{\pi_1(X,x)})$  holds for every $\pi_1(X_0,x)$-module $M$ which is of one of the following forms
\begin{itemize}[leftmargin=* ,parsep=0cm,itemsep=0cm,topsep=0cm]
 \item (\ref{Th:TensorInvariantsHeart}.1) a torsion-free quotient of $ \SH^*(Y_{x},\Z_\ell)$;
 \item (\ref{Th:TensorInvariantsHeart}.2) a submodule of $ \SH^*(Y_{x},\Z_\ell)$ with torsion-free cokernel.
\end{itemize}
 \end{theoreme}
 
\section{Proof of Theorem \ref{Th:TensorInvariantsHeart}}
 \subsection{Proof of (\ref{Th:TensorInvariantsHeart}.0)}\label{Sec:TensorInvariants}
  We may assume that $x$  is a geometric point on $X$ over $x_0\in X_0(k_0)$  and  that $R^*f_*\Z_\ell$ is torsion-free (Fact \ref{Gabber}). In particular, we have the short exact sequence $$0\rightarrow R^*f_*\Z_\ell\stackrel{\ell}{\rightarrow}R^*f_*\Z_\ell \rightarrow  R^*f_*\F_\ell\rightarrow 0$$
and the assertion of Theorem \ref{Th:TensorInvariantsHeart} is equivalent to the fact that the canonical (injective) morphism 
$$\SH^0(X,R^*f_*\Z_\ell)\otimes \F_\ell\rightarrow \SH^0(X,R^*f_*\F_\ell) $$
is an isomorphism. This, in turn, amounts to showing that $\SH^1(X,R^*f_*\Z_\ell)[\ell]=0$. To show this, we compute - in two ways - the characteristic polynomial of $F(:=F_{k_0})$ acting on $\SH^1(X,R^wf_*\Z_\ell)[\ell]$. \\
\indent On the one hand, we have $$\SH^1(X,R^wf_*\Z_\ell)[\ell]\twoheadleftarrow \SH^0(X,R^wf_*\F_\ell)\simeq \SH^w(Y_x, \F_\ell)^{\pi_1(X,x)}\subset \SH^w(Y_x,\F_\ell)\twoheadleftarrow \SH^w(Y_x,\Z_\ell)\hookrightarrow \SH^w(Y_x,\Q_\ell),$$
which shows that $F$ acting on $\SH^1(X,R^wf_*\Z_\ell)[\ell]$ is killed by  
 $$P_w:=P_{x_0}^{R^wf_{0*}\Q_\ell}=\Sdet(TId-F,\SH^w(Y_x,\Q_\ell))\in \Q[T],$$
which is independent of $\ell(\not=p)$ and whose roots are $|k_0|$-Weil numbers of weight $w$  \cite[Cor. 3.3.9]{Weil2}. Here we use that $x_0\in X_0(k_0)$ hence (Subsection \ref{Subsec:Not2}), that the action of  $F$  on $\SH^0(X,R^wf_*\F_\ell)$  identifies with the restriction of the action of $F_{x_0}$ on $ \SH^w(Y_x,\F_\ell)$.\\
\indent On the other hand, from Lemma \ref{Lem:Tech} below, the characteristic polynomial of $F$ acting on $\SH^1(X,R^wf_*\Z_\ell)[\ell] $ divides the characteristic polynomial of $F$ acting on  $\SH^1(X,R^wf_*\Z_\ell)\otimes\F_\ell$. As we also have a canonical $F$-equivariant embedding   $\SH^1(X,R^wf_*\Z_\ell)\otimes\F_\ell\subset \SH^1(X,R^wf_*\F_\ell)$, Lemma \ref{Lem:Weights} below shows that there exists $P_{\geq w+1}\in \Q[T]$, which is independent of $\ell(\not=p)$, whose roots are $|k_0|$-Weil numbers of weights $\geq w+1$ and such that $F$ acting on $\SH^1(X,R^wf_*\Z_\ell)[\ell]$ is killed by $P_{\geq w+1}$ for $\ell\gg 0$.\\
\indent The conclusion thus follows from the fact that $P_w,P_{\geq w+1}\in \Q[T]$ are coprime hence that $\overline{P}^\ell_w,\overline{P}^\ell_{\geq w+1}\in \F_\ell[T]$ are coprime as well for $\ell\gg 0$. $\square$\\

\begin{lemme}\label{Lem:Tech} Let $H$ be a finitely generated $\Z_\ell$-module equipped with a $\Z_\ell$-linear automorphism $F$.  Then the characteristic polynomial of $F$ acting on $H[\ell]$ always divides  the characteristic polynomial of $F$ acting on $H\otimes \F_\ell$.
\end{lemme}
\noindent\textit{Proof.}  As $\Z_\ell$ is a P.I.D., the short exact sequence of $\Z_\ell$-modules $0\rightarrow H_{tors}\rightarrow H\rightarrow H/H_{tors}\rightarrow 0$ always splits. In particular 
$H_{tors}\otimes \F_\ell\subset H\otimes \F_\ell$. So it is enough to show that the characteristic polynomial of $F$ acting on $H[\ell]$ always divides  the characteristic polynomial of $F$ acting on $H_{tors}\otimes \F_\ell$. Hence we may assume $H$ is finite. Then  the exact sequence of $\Z_\ell[F]$-modules
$$0\rightarrow H[\ell]\rightarrow H\stackrel{\ell}{\rightarrow}H\rightarrow H\otimes\F_\ell\rightarrow 0$$
shows that, in the Grothendieck group of $\Z_\ell[F]$-modules of finite length we have $[H[\ell]]=[H\otimes\F_\ell]$. In particular, $H[\ell]$ and $H\otimes\F_\ell$ have the same $F$-semisimplification. $\square$\\

\begin{lemme}\label{Lem:deJong} Let $\Lambda_\ell$ denote   $\Q_\ell$ (resp. $\F_\ell$). Let $Y_0$ be a scheme, separated and of finite type over a finite field $k_0$ and let $Y $ denote the base-change of $Y_0$ over $k:=\overline{k}_0$. Then for every integer $w\geq 0$ there exists $P_{\leq w,Y_0}\in \Q[T]$, which is independent of $\ell(\not=p)$ and whose roots are $|k_0|$-Weil numbers of weights $\leq w$ and such that $F$ acting on $\SH_c^w(Y,\Lambda_\ell)$ is killed by $P_{\leq w,Y_0}\in\Q[T]$ (resp. for $\ell\gg 0$ (depending on $Y$)).
\end{lemme}

\noindent\textit{Proof.} (See  \cite[Lemma 4.1]{Geisser}). The assertion holds for smooth proper $Y$ (\cite[Cor. 3.3.9]{Weil2} for $\Lambda_\ell=\Q_\ell$ plus Fact \ref{Gabber} for $\Lambda_\ell=\F_\ell$). Also, if the assertion holds over a  finite extension $k_0\hookrightarrow k_0'$ then it holds over $k_0$. Indeed, if $d:=[k_0':k_0]$ and $F_{k_0'}=F_{k_0}^d$ acting on $\SH_c^w(Y,\Lambda_\ell)$ is killed by $P_{\leq w,Y_0\otimes_{k_0} k_0'}\in\Q[T]$  whose roots are $|k_0'|$-Weil numbers of weights $\leq w$ then $F_{k_0}$ is killed by $P_{\leq w,Y_0}:=P_{\leq w,Y_0\otimes_{k_0} k_0'}(T^d)\in\Q[T]$ whose roots are $|k_0|$-Weil numbers of weights $\leq w$. Thus, in the following, we will implicitly allow finite field extensions of the base field (for instance the divisor $D$, alteration $Y'$ \textit{etc.} introduced below may only be defined over a finite extension of $k_0$, but this does not affect the argument). Eventually, by topological invariance of \'etale cohomology, we may assume that $Y$ is reduced.\\
\indent We proceed by  induction on the dimension of $Y$. The $0$-dimensional case is straightforward. Assume the assertion of the lemma holds for $\leq r$-dimensional reduced schemes, separated and of finite type over $k_0$. Fix an  $(r+1)$-dimensional scheme $Y_0$, separated and of finite type over $k_0$.  Write $D$ for the union of the intersections of the pairs of distinct irreducible components of $Y$. Then the localization exact sequence for cohomology with compact support
$$\cdots\rightarrow \SH^w_c(Y\setminus D,\Lambda_\ell)\rightarrow \SH^w_c(Y,\Lambda_\ell)\rightarrow\SH^w_c(D,\Lambda_\ell)\rightarrow\cdots$$
and the induction hypothesis for $D$ show that, without loss of generality, $Y$ may be assumed to be integral. 
From de Jong's alterations theorem \cite[Thm. 4.1]{dJ}, there exists a generically \'etale alteration $\phi:Y'\rightarrow Y$ and an open embedding $Y'\hookrightarrow Y^{'cpt}$ into a scheme $Y^{'cpt}$ smooth and projective over $k$ such that $Y^{'cpt}\setminus Y'$ is a strict normal crossing divisor. Fix a non-empty open subscheme $\emptyset\not= U\hookrightarrow Y$ such that $U':=Y'\times_YU\rightarrow U$ is finite \'etale and write $D:=Y\setminus U$. Again, the localization exact sequence for cohomology with compact support and the induction hypothesis for $D$ show that it is enough to prove the claim for $U$. As $U'\rightarrow U$ is a finite \'etale morphism of degree say $\delta$, one has the trace morphism which induces the multiplication-by-$\delta$ morphism
$$\delta\cdot :\SH_c^w(U,\Lambda_\ell)\stackrel{res}{\rightarrow}\SH_c^w(U',\Lambda_\ell)\stackrel{tr}{\rightarrow}\SH_c^w(U,\Lambda_\ell).$$
So, as soon as $\ell >\delta $, $\SH_c^w(U,\Lambda_\ell)$ is a direct factor of $\SH^w(U',\Lambda_\ell)$ as an $F$-module. Hence it is enough to prove the claim for $U'$. Write $D':=Y^{' cpt}\setminus U'$. Then the  localization exact sequence for cohomology with compact support
$$\cdots\rightarrow \SH^{w-1}_c(D',\Lambda_\ell)\rightarrow \SH^w_c(U',\Lambda_\ell)\rightarrow\SH^w_c(Y^{' cpt},\Lambda_\ell)\rightarrow\cdots,$$
the induction hypothesis for $D'$ and the fact that the assertion of the lemma holds for the smooth projective scheme $Y^{' cpt}$ yield the conclusion. $\square$

 \begin{lemme}\label{Lem:Weights} Let $\Lambda_\ell$ denote $\Q_\ell$ (resp. $\F_\ell$). With the notation of Theorem \ref{Th:TensorInvariantsHeart} there exists $P_{\geq w+1}\in\Q[T]$ whose roots are $|k_0|$-Weil numbers of weights $\geq w+1$ and such that $F$ acting on $\SH^1(X,R^wf_*\Lambda_\ell)$ is killed by $P_{\geq w+1}$ (resp. for $\ell\gg 0$). 
 \end{lemme}
  
  \noindent \textit{Proof.} One may assume that $Y_0$ is connected, hence irreducible. Then $Y$ is equidimensional, say, of dimension $d_Y$. From the Leray spectral sequence $E_2^{v,w}=R^vg_*R^wf_*\Lambda_\ell\Rightarrow R^{v+w}(gf)_*\Lambda_\ell$ for 
$Y_0\stackrel{f_0}{\rightarrow}X_0\stackrel{g_0}{\rightarrow}\spec(k_0)$, one sees that $\SH^1(X ,R^wf_*\Lambda_\ell)=E_2^{1,w}=E_\infty^{1,w} $ (recall that $X$ is a curve) is a subquotient of $R^{1+w}(gf)_*\Lambda_\ell=\SH^{1+w}(Y,\Lambda_\ell)\simeq \SH^{2d_Y-w-1}_c(Y,\Lambda_\ell)(d_Y)^\vee$ (the second isomorphism is Poincar\'e duality). Now, take $P_{\leq 2d_Y-w-1, Y_0}\in\Q[T]$ as in Lemma \ref{Lem:deJong} and let $\delta$ denote its degree. Then, for $\Lambda_\ell=\Q_\ell$ (resp. $\F_\ell$), $F$ acting on $\SH^1(X ,R^wf_*\Lambda_\ell)$ is killed by $T^\delta P_{\leq 2d_Y-w-1, Y_0}(|k_0|^{d_Y} T^{-1})\in\Q[T]$, whose roots are $|k_0|$-Weil numbers of weights $\geq w+1$, as desired. $\square$\\

\subsection{Proof of  (\ref{Th:TensorInvariantsHeart}.1)} We may assume $X_0$ is affine. From (\ref{Th:TensorInvariantsHeart}.0) we have an $F$-equivariant injective morphism $$ H^1(\pi_1(X,x),H^w(Y_x,\Z_\ell))\hookrightarrow H^1(\pi_1(X,x),H^w(Y_x,\Z_\ell))\otimes\Q_\ell=H^1(X,R^wf_*\Q_\ell).$$
From the  $\Lambda_\ell=\Q_\ell$ case of Lemme \ref{Lem:Weights}, $F$ acting on $H^1(X,R^wf_*\Q_\ell)$ is killed by $P_{\geq w+1}\in \Q[T]$ independent of $\ell$. So the same is true for $H^1(\pi_1(X,x),H^w(Y_x,\Z_\ell))$. On the other hand, as $X$ has cohomological dimension $\leq 1$, the canonical $F$-equivariant morphism
$$H^1(\pi_1(X,x),H^w(Y_x,\Z_\ell))\rightarrow H^1(\pi_1(X,x),M)$$
is surjective. This shows that\\

\noindent (5.2.1) $F$ acting on $H^1(\pi_1(X,x),M)$ is killed by a polynomial $P_{\geq w+1}\in \Q[T]$ independent of $\ell$ and $M$ and whose roots are $|k_0|$-Weil numbers of weights $\geq w+1$ .\\

\noindent The assertion (5.2.1) implies that $H^1(\pi_1(X,x),M)[\ell]$ is also killed by $P_{\geq w+1}$. To conclude, consider the diagram
$$\xymatrix{0\ar[r]&M^{\pi_1(X,x)}\otimes\F_\ell\ar[r] &(M\otimes\F_\ell)^{\pi_1(X,x)}\ar[r]\ar@{_{(}->}[d]&H^1(\pi_1(X,x),M)[\ell]\ar[r]&0\\
H^w(Y_x,\Z_\ell)\ar@{->>}[r]&M\ar@{->>}[r]&M\otimes\F_\ell,}$$
which shows that $F$ acting on $(M\otimes\F_\ell)^{\pi_1(X,x)}$ is killed by a polynomial $P_{w}\in\Q[T]$ independent of $\ell$ and $M$ and whose roots are $|k_0|$-Weil numbers of weight $w$.\\

\subsection{Proof of  (\ref{Th:TensorInvariantsHeart}.2)}\label{ProofTI2} Let $\Pi_{\ell^\infty}$ denote  the image of $\pi_1(X,x)$ acting on $H^w(Y_x,\Z_\ell)$. From Fact \ref{DeligneLP}, $\Pi_{\ell^\infty}^{ab}$ is finite (see \textit{e.g.} \cite[Thm. 5.7]{UOI1}). In particular, $\Pi_{\ell^\infty}$ acts on $\det(M)$ through a finite quotient, which has to be of order dividing $\ell-1$. But, on the other hand $\Pi_{\ell^\infty}$ is generated by its $\ell$-Sylow subgroups (Fact \ref{JNTI}) so $\Pi_{\ell^\infty}$ acts trivially on $\det(M)$.  This shows that  the canonical isomorphism 
$$M\tilde{\rightarrow}\SHom(\Lambda^{m-1}M,\hbox{\rm det}(M))=(\Lambda^{m-1}M)^\vee\otimes \hbox{\rm det}(M),\;\; v\mapsto v\wedge-$$
(where $m$ denotes the $\Z_\ell$-rank of $M$) induces a $\pi_1(X,x)$-equivariant isomorphism $M\simeq (\Lambda^{m-1}M)^\vee$.  As $M\hookrightarrow H^w(Y_x,\Z_\ell)$ has torsion-free cokernel, $M^\vee$ is a torsion-free $\pi_1(X_0,x)$-quotient of $H^w(Y_x,\Z_\ell)^\vee$ hence $ \Lambda^{m-1}M^\vee$ is a torsion-free $\pi_1(X_0,x)$-quotient of $\Lambda^{m-1}H^w(Y_x,\Z_\ell)^\vee $, which is itself a $\pi_1(X_0,x)$-quotient of $\SH^{(m-1)(2d_f-w)}(Y_x^{[m-1]},\Z_\ell)((m-1)d_f)$ for $\ell\gg 0$. So we are reduced to (\ref{Th:TensorInvariantsHeart}.1).  $\square$\\

 \section{Summary}\label{Summary}  
 
 Our goal in the remaining parts of this paper  is to prove Theorem \ref{Th:GSSMax} and Theorem \ref{Th:GSS}. We fix the notation and conventions which will be used from now on and review the information we have collected so far.\\
 
 \noindent Let $k_0$ be a field finitely generated over $\F_p$ and contained in an algebraically closed field $k$ and let $f_0:Y_0\rightarrow X_0$ be a smooth proper morphism of $k_0$-schemes, with $X_0$ smooth separated, geometrically connected over $k_0$. Let $f:Y\rightarrow X$ denote the base-change of $f_0:Y_0\rightarrow X_0$ over $k$.    Assume $\ell\gg 0$ so that $R^*f_*\Z_\ell$ is torsion-free of constant rank $r$ (Fact \ref{Gabber}).   \\

\noindent 
\noindent Let $\Pi_{\ell^\infty}$ (resp. $\Pi_{0 \ell^\infty}$) denote the image of $\pi_1(X,x)$ (resp.  $\pi_1(X_0,x)$) acting on $\SH^*(Y_x,\Z_\ell)=:H_{\ell^\infty}$ and let $\Pi_{\ell}$ denote the image of $\pi_1(X,x)$ acting on $\SH^*(Y_x,\F_\ell)=:H_\ell$. Write $V_{\ell^\infty}:=H_{\ell^\infty}\otimes_{\Z_\ell}\Q_\ell$.  \\

\noindent  Let   $\frak{G}_{\ell^\infty}\hookrightarrow \SGL_{H_{\ell^\infty}}$ denote the Zariski closure of $\Pi_{\ell^\infty}$ (endowed with the reduced subscheme structure), $\mathcal{G}_{\ell^\infty}:=\frak{G}_{\ell^\infty,\Q_\ell}\hookrightarrow \SGL_{V_{\ell^\infty}}$ its generic fiber, $\mathcal{G}_\ell:=\frak{G}_{\ell^\infty,\F_\ell}\hookrightarrow \SGL_{H_\ell}$ its special fiber. The scheme $\mathcal{G}_{\ell^\infty}$ coincides with the Zariski closure of $\Pi_{\ell^\infty}$ in $\SGL_{V_{\ell^\infty}}$. Let also  $\mathcal{G}_{0\ell^\infty}$ denote the Zariski closure of $\Pi_{0\ell^\infty}$ in $\SGL_{V_{\ell^\infty}}$.  By construction\\

\begin{itemize}[leftmargin=* ,parsep=0cm,itemsep=0cm,topsep=0cm]
 \item (\ref{Summary}.1) $\frak{G}_{\ell^\infty}\hookrightarrow \SGL_{H_{\ell^\infty}}$ is  flat  over $\Z_\ell$. In particular, $\Sdim(\mathcal{G}_{\ell})=\Sdim(\mathcal{G}_{\ell^\infty})$.\\  
\end{itemize}

\noindent  From Lemma \ref{Lem:Localization},   we may assume (Fact \ref{DeligneLP}, Fact \ref{JNTI})\\

\begin{itemize}[leftmargin=* ,parsep=0cm,itemsep=0cm,topsep=0cm]
 \item (\ref{Summary}.2.1) $\Pi_\ell $ is perfect and generated by its order $\ell$ elements for  $\ell\gg 0$ ;
  \item (\ref{Summary}.2.2)  $\mathcal{G}_{\ell^\infty}$ is connected semisimple and $\mathcal{G}_{0\ell^\infty}$ is connected for every $\ell\not=p$.\\  
\end{itemize}
 
 \noindent Eventually, Theorem \ref{Th:TensorInvariants} reads\\

\begin{itemize}[leftmargin=* ,parsep=0cm,itemsep=0cm,topsep=0cm]
 \item (\ref{Summary}.3)  For every integers $d,d^\vee\geq0$, partitions $\lambda, \lambda^\vee$ of $d, d^\vee$   and $\ell\gg 0$ (depending on  $d, d^\vee$) the property   $(Inv,M)$  holds for 
\begin{itemize}
\item (\ref{Summary}.3.1) Every $\Pi_{0\ell^\infty}$-module  quotient $S_{\lambda,\lambda^\vee}(H_{\ell^\infty})\twoheadrightarrow M$ which is torsion-free; 
\item  (\ref{Summary}.3.2) Every  $\Pi_{0\ell^\infty}$-submodule $M\hookrightarrow S_{\lambda,\lambda^\vee}(H_{\ell^\infty})$ with torsion-free cokernel.\\
\end{itemize}  
 \end{itemize}

\begin{center}  \textbf{ \fontfamily{lmr}\selectfont{PART II: SEMISIMPLICITY VERSUS MAXIMALITY}} \end{center}

 \section{Structure of $\mathcal{G}_\ell$; first reformulations of (\ref{Th:GSS})}
 
 \subsection{Group-theoretical preliminaries} Let $\Lambda_\ell$ denote $\Z_\ell$ or $\F_\ell$. Given a closed subgroup $\Pi$ of $\SGL_r(\Lambda_\ell)$, write $\Pi^+\subset \Pi $ for the (normal closed) subgroup of $\Pi $ generated by its $\ell$-Sylow subgroups. Given a finitely generated $\Lambda_\ell$-module $H$, write $$H^\otimes:=\bigoplus_{s,t \geq 0}H^{\otimes s}\otimes (H^\vee)^{\otimes t}.$$
 For an integer $d\geq 1$, set  $$T^{\leq d}(H):=\bigoplus_{s,t\leq d}H^{\otimes s}\otimes (H^\vee)^{\otimes t}\subset H^\otimes.$$

\noindent Let $H_\ell$ be an $r$-dimensional $\F_\ell$-vector space. Given a subgroup $\Pi_\ell\subset \SGL(H_\ell)$, let $\widetilde{\Pi}_\ell\hookrightarrow \SGL_{H_\ell}$ denote its algebraic envelope, in the sense of Nori  \cite{Nori} that is the algebraic subgroup generated by the one-parameter groups
$$\begin{tabular}[t]{llll}
$\phi_g$:&$\A_{\F_\ell}^1$&$\rightarrow$&$ \SGL_{H_\ell}$\\
&$t$&$\mapsto$&$\exp(t\log(g))$ 
\end{tabular}$$
for $g\in\Pi_\ell $  of order $\ell$. Here $\exp(n):=\sum_{0\leq i\leq \ell-1}\frac{n^i}{i!}$ for a nilpotent $n\in \SEnd(H_\ell)$ and $\log(u):=-\sum_{1\leq i\leq \ell-1}\frac{(1-u)^i}{i}$ for a unipotent $u\in \SGL(H_\ell)$.\\

\noindent By construction $\widetilde{\Pi}_\ell$ is a smooth algebraic subgroup of $ \SGL_{H_\ell}$, connected and generated by its unipotent  subgroups. Furthermore, the following hold.
\begin{lemme}\label{Fact:NoriEnvelope} \textit{}
\begin{itemize}[leftmargin=* ,parsep=0cm,itemsep=0cm,topsep=0cm] 
\item (\ref{Fact:NoriEnvelope}.1) For $\ell\gg 0$ depending only on $r$, we have $\widetilde{\Pi}_\ell(\F_\ell)^+=\Pi_\ell^+$. In particular,  $\Pi_\ell^+$-submodules and $\widetilde{\Pi}_\ell$-submodules of $H_\ell$ coincide. 
\item (\ref{Fact:NoriEnvelope}.2) There exists  $d\geq 1$ depending only on $r$ such that $\widetilde{\Pi}_\ell$ is  the stabilizer of $(T^{\leq d}(H_\ell))^{\widetilde{\Pi}_\ell}$ in $\SGL_{H_\ell}$, for $\ell\gg 0$ depending only on $r$.
\item (\ref{Fact:NoriEnvelope}.3) For every $d\geq 1$ and $\ell\gg 0$ depending only on $d$, $r$ we have $$(T^{\leq d}(H_\ell))^{\widetilde{\Pi}_\ell}=(T^{\leq d}(H_\ell))^{\Pi_\ell^+}.$$
\end{itemize}
\end{lemme}
\noindent\textit{Proof.} The first part of (\ref{Fact:NoriEnvelope}.1) is \cite[Thm. B]{Nori} and the second part follows from the first part by construction of $\widetilde{\Pi}_\ell$. The assertion (\ref{Fact:NoriEnvelope}.2) follows from \cite[Lem. 4.1]{TI}. For (\ref{Fact:NoriEnvelope}.3), the inclusion $(T^{\leq d}(H_\ell))^{\widetilde{\Pi}_\ell}\subset (T^{\leq d}(H_\ell))^{\Pi_\ell^+}$  holds as soon as $\ell\geq r$ (recall that for $\ell\geq r$ the only elements in $\SGL_r(\F_\ell)$ of order a power of $\ell$ are those of order $\ell$). For the opposite inclusion, fix an isomorphism $H_\ell\tilde{\rightarrow}\F_\ell^{\oplus r}$. Then for every $v\in (T^{\leq d}(H_\ell))^{\Pi_\ell^+}$ and $g\in \Pi_\ell$ of order $\ell$ each component of the vector  equation $$\exp(t\log(g))v-v=0$$
is a polynomial in $t$ with degree $\leq 2d(r-1)$ and has at least $\ell$ distinct roots. So for $\ell> 2d(r-1)$ the image of $\phi_g$ is contained in the stabilizer of $v$ in $\SGL_{H_\ell}$. This shows $(T^{\leq d}(H_\ell))^{\widetilde{\Pi}_\ell}\supset (T^{\leq d}(H_\ell))^{\Pi_\ell^+}$. $\square$\\

\noindent A semisimple group scheme over $\Z_\ell$ is a smooth affine group scheme whose geometric fibers are  connected semisimple algebraic groups.  Then all the geometric fibers have the same root data \cite[XXII, Prop. 2.8]{SGA3}. We say that a semisimple group scheme over $\Z_\ell$ is simply connected if its fibers are. Furthermore, we have

\begin{lemme}\label{Fact:Frattini} Let $\frak{G}$ be a simply connected semisimple group scheme over $\Z_\ell$. Then  $\frak{G}(\Z_\ell)=\frak{G}(\Z_\ell)^+$.
 
\end{lemme} 
\noindent\textit{Proof.} This follows from the fact that the kernel of the reduction-modulo-$\ell$ morphism $\frak{G}(\Z_\ell)\rightarrow \frak{G}(\F_\ell)$ is a pro-$\ell$ group and that, as $\frak{G}_{\F_\ell}$ is simply connected, $\frak{G}(\F_\ell)=\frak{G}(\F_\ell)^+$  \cite[\S 12]{Steinberg}. $\square$

\subsection{Structure of $\mathcal{G}_\ell$ and weak maximality} We are now able to  prove the main result of this section.    
\begin{theoreme}\label{Cor:SSConsequence1} For $\ell\gg 0$, we have $\widetilde{\Pi}_\ell=\mathcal{G}_{\ell}$. In particular,
\begin{itemize}[leftmargin=* ,parsep=0cm,itemsep=0cm,topsep=0cm]
\item (\ref{Cor:SSConsequence1}.1) For $\ell\gg 0$, $\mathcal{G}_{\ell}(\F_\ell)^+$ is perfect and $\mathcal{G}_{\ell}=\widetilde{\mathcal{G}_{\ell}(\F_\ell)}$;
\item (\ref{Cor:SSConsequence1}.2) (Weak maximality) There exists an integer $C_r\geq 1$ depending only on $r$ such that, for $\ell\gg 0$, $[\frak{G}_{\ell^\infty}(\Z_\ell):\Pi_{\ell^\infty}]\leq C_r$ (hence $\Pi_{\ell^\infty}=\frak{G}_{\ell^\infty}(\Z_\ell)^+$).
\end{itemize}
\end{theoreme}
\noindent\textit{Proof.}  Assume $\ell\gg 0$ so that  (\ref{Summary}.2.1), (\ref{Summary}.2.2) and the conclusions of  Lemma \ref{Fact:NoriEnvelope}  hold. By construction, $(H_{\ell^\infty}^\otimes)^{\Pi_{\ell^\infty}}=(H_{\ell^\infty}^\otimes)^{\frak{G}_{\ell^\infty}}$  and $\frak{G}_{\ell^\infty}$ is contained in the stabilizer of $(H_{\ell^\infty}^\otimes)^{\Pi_{\ell^\infty}}$.  As stabilizers commute with arbitrary base-changes, $\mathcal{G}_\ell$ is  contained in the stabilizer of $(H_{\ell^\infty}^\otimes)^{\Pi_{\ell^\infty}}\otimes\F_\ell$  hence, \textit{a fortiori} in the stabilizer of $T^{\leq d}(H_{\ell^\infty})^{\Pi_{\ell^\infty}}\otimes\F_\ell$ for every integer $d\geq 1$.  On the other hand, from (\ref{Fact:NoriEnvelope}.2), there exists an integer $d\geq 1$ depending only on $r$ such that $\widetilde{\Pi}_\ell$ is the stabilizer of $(T^{\leq d}(H_\ell))^{\widetilde{\Pi}_\ell}$ in $\SGL_{H_\ell}$. Then,  up to increasing $\ell$,  we have  $T^{\leq d}(H_{\ell^\infty})^{\Pi_{\ell^\infty}}\otimes\F_\ell=T^{\leq d}(H_\ell)^{\Pi_\ell}$ (Theorem \ref{Th:TensorInvariants}). By (\ref{Summary}.2.1) and  (\ref{Fact:NoriEnvelope}.3) this shows that $\mathcal{G}_\ell\subset \widetilde{\Pi}_\ell$. For the opposite inclusion, as $\tilde{\Pi}_\ell$ is integral (being smooth and connected), it is enough to show that $\Sdim(\mathcal{G}_{\ell})\geq \Sdim(\widetilde{\Pi}_\ell)$. This follows from $\Sdim(\mathcal{G}_{\ell^\infty})\geq \Sdim(\widetilde{\Pi}_\ell)$  \cite[Thm. 7]{LarsenpAdicNori} and  (\ref{Summary}.1). Then (\ref{Cor:SSConsequence1}.1) follows from (\ref{Summary}.2.1) and (\ref{Fact:NoriEnvelope}.1) while the first part of (\ref{Cor:SSConsequence1}.2) follows from \cite[Thm. 7 (3)]{LarsenpAdicNori}. The assertion $\Pi_{\ell^\infty}=\frak{G}_{\ell^\infty}(\Z_\ell)^+$ then follows from the fact  that $\Pi_{\ell^\infty}^+\subset\frak{G}_{\ell^\infty}(\Z_\ell)^+\subset \Pi_{\ell^\infty} $  for  $\ell>C_r$ and (\ref{Summary}.2.1).
$\square$
 
 \begin{corollaire}\label{Lem:SmoothRS} The subgroup scheme $\frak{G}_{\ell^\infty}\hookrightarrow \SGL_{H_{\ell^\infty}}$ is connected, smooth  over $\Z_\ell$ for $\ell\gg 0$.
\end{corollaire}
 \noindent\textit{Proof.} The key point is that $\widetilde{\Pi}_\ell=\mathcal{G}_\ell$ for $\ell\gg 0$ (Theorem \ref{Cor:SSConsequence1}). Then, the assertion about smoothness follows from the fact that $\frak{G}_{\ell^\infty}$ is flat (\ref{Summary}.1) and of finite type over $\Z_\ell$  and from the smoothness of $\widetilde{\Pi}_\ell$ (as observed in the paragraph before Lemma \ref{Fact:NoriEnvelope}) while the assertion about connectedness follows from (\ref{Summary}.2.2) and the connectedness of  $\widetilde{\Pi}_\ell$. $\square$

  \subsection{First reformulations of (\ref{Th:GSS})} Using Corollary \ref{Lem:SmoothRS}, we obtain the following reformulations of (\ref{Th:GSS}).
 
\begin{corollaire}\label{Cor:SSConsequence2}  The following assertions are equivalent:
\begin{itemize}[leftmargin=* ,parsep=0cm,itemsep=0cm,topsep=0cm]
\item (\ref{Cor:SSConsequence2}.1) the action  of $\Pi_\ell$ on $H_\ell$ is semisimple for $\ell\gg 0$;
\item (\ref{Cor:SSConsequence2}.2) the action of $\widetilde{\Pi}_\ell$ (equivalently $\mathcal{G}_\ell$) on $H_\ell$ is semisimple for $\ell\gg 0$;
\item (\ref{Cor:SSConsequence2}.3) $\widetilde{\Pi}_\ell$ (equivalently $\mathcal{G}_\ell$) is semisimple for  $\ell\gg 0$;
\item (\ref{Cor:SSConsequence2}.4) $\frak{G}_{\ell^\infty}$ is a semisimple group scheme over $\mathbb{Z}_\ell$ for  $\ell\gg 0$.
\end{itemize}
\end{corollaire}
\noindent\textit{Proof.}  Assume $\ell\gg 0$ so that  (\ref{Summary}.2.1), (\ref{Summary}.2.2), (\ref{Summary}.3) and the conclusions of  Lemma \ref{Fact:NoriEnvelope} and  Corollary \ref{Lem:SmoothRS} hold. The equivalence (\ref{Cor:SSConsequence2}.1) $\Leftrightarrow$ (\ref{Cor:SSConsequence2}.2) follows from (\ref{Summary}.2.1) and (\ref{Fact:NoriEnvelope}.1). The fact that  (\ref{Cor:SSConsequence2}.2) implies that $\tilde{\Pi}_\ell$ is reductive is standard.  $\tilde{\Pi}_\ell$ is then automatically semisimple since it it generated by its unipotent subgroups. This shows
 (\ref{Cor:SSConsequence2}.2) $\Rightarrow$ (\ref{Cor:SSConsequence2}.3). The equivalence (\ref{Cor:SSConsequence2}.3) $\Leftrightarrow$ (\ref{Cor:SSConsequence2}.4) is by definition (since (\ref{Summary}.2.2), Corollary \ref{Lem:SmoothRS} holds).
The  implication (\ref{Cor:SSConsequence2}.3) $\Rightarrow$ (\ref{Cor:SSConsequence2}.2)  follows for instance from   \cite[Prop. 3.2]{Jantzen} (see also \cite[Thm. 3.5]{LarsenSS}).
 $\square$

\section{Semisimplicity versus maximality}  Let  $\mathcal{G}$ be a connected semisimple group over $\Q_\ell$. Write $p^{sc}:\mathcal{G}^{sc}\rightarrow \mathcal{G}$ and $p^{ad}:\mathcal{G} \rightarrow \mathcal{G}^{ad}$ for the simply connected cover and adjoint quotient of $\mathcal{G}$ respectively. Recall that the  Bruhat-Tits building \cite{Tits}  $\mathcal{B}:=\mathcal{B}(\mathcal{G}^{sc},\Q_\ell)$ is equipped with a natural action of   $\mathcal{G}^{ad}(\Q_\ell)$ and that   $\mathcal{G}(\Q_\ell)$ acts on $\mathcal{B}$ through its image in $\mathcal{G}^{ad}(\Q_\ell)$.  There is a bijective correspondence between
\begin{itemize} 
\item semisimple models $\frak{G}$ of $\mathcal{G}$ over $\Z_\ell$;
\item hyperspecial points  $b\in \mathcal{B}$,
\end{itemize}
given by $\frak{G}\rightarrow \mathcal{B}^{\frak{G}(\Z_\ell)}$ \cite[3.8.1]{Tits}.
Also given an isogeny $\phi:\mathcal{G} \rightarrow \mathcal{G}' $ and if $\frak{G}_b$, $\frak{G}'_b$ are respectively the semisimple models over $\Z_\ell$ of $ \mathcal{G}$ and $\mathcal{G}'$ corresponding to a hyperspecial point $b\in \mathcal{B}$ then $\phi:\mathcal{G} \rightarrow \mathcal{G}'$ extends uniquely to a morphism $\phi_b:\frak{G}_b\rightarrow \frak{G}'_b$ of group schemes over $\Z_\ell$.\\

\noindent  A compact subgroup $\Pi\subset \mathcal{G}(\Q_\ell)$ of the form $\Pi=\frak{G}(\Z_\ell)$ for some  semisimple model $\frak{G}$ of $\mathcal{G}$ over $\Z_\ell$ (or, equivalently, such that $\Pi\subset \mathcal{G}(\Q_\ell)$ is the stabilizer of a hyperspecial point in $\mathcal{B}$) is called hyperspecial. Hyperspecial subgroups, when they exist, are the  compact subgroups of $\mathcal{G}(\Q_\ell)$ of maximal volume \cite[3.8.2]{Tits}. In particular, a $\mathcal{G}(\Q_\ell)$-conjugate of a hyperspecial subgroup is again hyperspecial.\\

\noindent  A compact subgroup $\Pi\subset \mathcal{G}(\Q_\ell)$ is called almost hyperspecial if $(p^{sc})^{-1}(\Pi)\subset \mathcal{G}^{sc}(\Q_\ell)$ is hyperspecial.

\begin{lemme}\label{Lem:Compare} Let $\mathcal{G}$ be a connected semisimple group over $\Q_\ell$.
Let  $\frak{G}$ (resp. $\frak{G}^{sc}$)  be a smooth, connected group scheme over $\Z_\ell$ with generic fiber $\mathcal{G}$ (resp. $\mathcal{G}^{sc}$). Assume 
\begin{itemize}[leftmargin=* ,parsep=0cm,itemsep=0cm,topsep=0cm]
\item (\ref{Lem:Compare}.1) $\frak{G}^{sc}$ is semisimple over $\Z_\ell$;
\item (\ref{Lem:Compare}.2)  $ (p^{sc})^{-1}(\frak{G}(\Z_\ell)^+)$ is a normal subgroup of $ \frak{G}^{sc}(\Z_\ell)$
such that $ \frak{G}^{sc}(\Z_\ell)/  (p^{sc})^{-1}(\frak{G}(\Z_\ell)^+)$ is abelian.
\end{itemize} Then for $\ell\gg 0$ depending only on the dimension of $\mathcal{G}$,  $\frak{G}$ is semisimple over $\Z_\ell$.
\end{lemme}
\noindent\textit{Proof.} For a profinite group $\Pi$ which is an extension of a finite group by a pro-$\ell$ group let $N(\Pi)$ denote the product of the orders of the groups (counted with multiplicities) appearing in the non-abelian part of the composition series of $\Pi$. Note that if $\Pi'\subset \Pi$ are two such groups then $N(\Pi')\leq N(\Pi)$.\\

\noindent Let $\frak{H}$ be a  connected, smooth, affine group scheme  over $\Z_\ell$; write $\mathcal{H}_\ell:=\frak{H}_{\F_\ell}$ for its special fiber. The following hold.\\

\begin{enumerate}[leftmargin=* ,parsep=0cm,itemsep=0cm,topsep=0cm]
\item  The non-abelian parts of the composition series of $\frak{H}(\Z_\ell)$  and $\frak{H}(\F_\ell)$ (resp.  $\frak{H}(\Z_\ell)^+$ and $\frak{H}(\F_\ell)^+$) coincide. This is because the reduction modulo-$\ell$ map $$\frak{H}(\Z_\ell)\rightarrow \frak{H}(\F_\ell)$$
is surjective with pro-$\ell$ kernel.\\
\item   Write   $\mathcal{H}_\ell^{ss}:=\mathcal{H}_\ell/R(\mathcal{H}_\ell)$, where $R(\mathcal{H}_\ell)$ is the solvable radical of $\mathcal{H}_\ell$. Then the non-abelian parts of the composition series of $\frak{H}(\F_\ell)$, $\mathcal{H}_\ell^{ss}(\F_\ell)$, $\mathcal{H}_\ell^{ss}(\F_\ell)^+$ and $\frak{H}(\F_\ell)^+$ coincide. Indeed, first, as  $\mathcal{H}_\ell^{ss}$ is a semisimple group, $\mathcal{H}_\ell^{ss}(\F_\ell)/\mathcal{H}_\ell^{ss}(\F_\ell)^+ $ is abelian for $\ell\gg 0$. More precisely, let  $\mu_\ell$ denote the kernel of the simply connected cover of $\mathcal{H}_\ell^{ss}$. Since   $\SH^1(\F_\ell,\mu_\ell(\overline{\F}_\ell))$ is of order dividing   $|\mu_\ell(\overline{\F}_\ell)|$ (\textit{e.g.} \cite[XIII, \S1, Prop. 1]{SerreCL}),  $\mathcal{H}_\ell^{ss}(\F_\ell)/\mathcal{H}_\ell^{ss}(\F_\ell)^+ $ embeds into  $\SH^1(\F_\ell,\mu_\ell(\overline{\F}_\ell))$ for $\ell$ prime to the order of $\mu_\ell$. The assertion then follows from the fact that  the rank of $\mathcal{H}_\ell^{ss}$ (hence the order of $\mu_\ell$) is bounded as $\ell$ varies. As a result,  the non-abelian parts of the composition series of $\mathcal{H}_\ell^{ss}(\F_\ell)$ and $\mathcal{H}_\ell^{ss}(\F_\ell)^+$ coincide. Next,  Lang's theorem \cite{Lang} gives a short exact sequence $$1\rightarrow  R(\mathcal{H}_\ell)(\F_\ell)\rightarrow \frak{H}(\F_\ell)\rightarrow \mathcal{H}_\ell^{ss}(\F_\ell)\rightarrow 1.$$
 Hence  the non-abelian parts of the composition series of $\frak{H}(\F_\ell)$ and $\mathcal{H}_\ell^{ss}(\F_\ell)$ coincide. Furthermore, the above short exact sequence
 induces a short exact sequence
$$1\rightarrow   R(\mathcal{H}_\ell)(\F_\ell)\cap\frak{H}(\F_\ell)^+\rightarrow \frak{H}(\F_\ell)^+\rightarrow \mathcal{H}_\ell^{ss}(\F_\ell)^+\rightarrow 1.$$
 Hence  the non-abelian parts of the composition series of $\frak{H}(\F_\ell)^+$ and $\mathcal{H}_\ell^{ss}(\F_\ell)^+$ coincide.  \\
 \item  Let $H_i$, $i\in I $ denote the almost simple factors of $\mathcal{H}_\ell^{ss} $. Then \cite[Main Thm.]{TitsSimple} the  non-abelian part of the composition series of  $\mathcal{H}_\ell^{ss}(\F_\ell) $ is precisely the family of the
$$ H_i(\F_\ell)^+/(Z(H_i(\F_\ell))\cap H_i(\F_\ell)^+),\; i\in I$$  
\noindent for $\ell\gg 0$. As the kernel and the cokernel of 
$\prod_{i\in I}H_i(\F_\ell)\rightarrow \mathcal{H}_\ell^{ss}(\F_\ell)$, the $ Z(H_i(\F_\ell)) $ and $H_i(\F_\ell)/H_i(\F_\ell)^+$, $i\in I$  all
have order bounded from above by a constant depending only on the rank of $\mathcal{H}_\ell^{ss} $, there exists a constant $c>0$ depending only on   the rank of $\mathcal{H}_\ell^{ss}$ such that 
$$ N(\mathcal{H}_\ell^{ss}(\F_\ell))\leq |\mathcal{H}_\ell^{ss}(\F_\ell)|\leq c N(\mathcal{H}_\ell^{ss}(\F_\ell)).$$
If $d$ denotes the dimension of $\mathcal{H}_\ell^{ss}$, this also implies \cite[Lemma 3.5]{Nori}
$$\frac{(\ell-1)^d}{c}\leq N(\mathcal{H}_\ell^{ss}(\F_\ell))\leq (\ell+1)^d .$$
 \end{enumerate}

\noindent As  $ \frak{G}^{sc}(\Z_\ell)/ (p^{sc})^{-1}(\frak{G}(\Z_\ell)^+)$ is abelian, the non-abelian parts of the composition series of  $ \frak{G}^{sc}(\Z_\ell)$ and $(p^{sc})^{-1}(\frak{G}(\Z_\ell)^+)$ coincide and as the kernel of $p^{sc}: (p^{sc})^{-1}(\frak{G}(\Z_\ell)^+)\rightarrow \frak{G}(\Z_\ell)^+$ is abelian, we have
$$N( \frak{G}^{sc}(\F_\ell))=N( \frak{G}^{sc}(\Z_\ell))=N((p^{sc})^{-1}(\frak{G}(\Z_\ell)^+))\leq N(\frak{G}(\Z_\ell)^+)=N(\frak{G}(\F_\ell)^+)=N(\frak{G}(\F_\ell))=N(\frak{G}_{\F_\ell}^{ss}(\F_\ell)).$$
\noindent Let $d$ denote the common dimension of $\frak{G}_{\F_\ell}$ and $\frak{G}^{sc}_{\F_\ell}$ and let $d^{ss}$ denote the dimension of $\frak{G}_{\F_\ell}^{ss}$. Then we have
$$ \frac{(\ell-1)^d}{c}\leq \frac{|\frak{G}^{sc}(\F_\ell)|}{c}  \leq N(\frak{G}^{sc}(\F_\ell))\leq  N(\frak{G}_{\F_\ell}^{ss}(\F_\ell))\leq  (\ell+1)^{d^{ss}} . $$
When $\ell\rightarrow +\infty$, this forces $d^{ss}=d$, as desired. $\square$\\

\begin{corollaire}\label{Cor:SSConsequence3} The assertions of Corollary \ref{Cor:SSConsequence2} are also equivalent to

\begin{itemize}[leftmargin=* ,parsep=0cm,itemsep=0cm,topsep=0cm]
\item (\ref{Cor:SSConsequence3}) $\Pi_{\ell^\infty}\subset \mathcal{G}_{\ell^\infty}(\Q_\ell)$ is an almost hyperspecial subgroup for $\ell\gg 0$.
\end{itemize} 
\end{corollaire}

\noindent\textit{Proof.}   Assume $\ell\gg 0$ so that (\ref{Summary}.2.1), (\ref{Summary}.2.2), (\ref{Summary}.3), (\ref{Cor:SSConsequence1}.2) and the conclusion of Corollary \ref{Lem:SmoothRS} hold.  Assume  (\ref{Cor:SSConsequence2}.4) holds. Then  $\frak{G}_{\ell^\infty}$ corresponds to a hyperspecial point $b\in \mathcal{B} $, which also gives rise to a simply connected semisimple model $\frak{G}_{\ell^\infty}^{sc}$ over $\Z_\ell$ of the simply connected cover $p^{sc}:\mathcal{G}_{\ell^\infty }^{sc}\rightarrow \mathcal{G}_{\ell^\infty }$ with the property that $\frak{G}_{\ell^\infty}^{sc}(\Z_\ell)$ is the stabilizer of $b$ in $\mathcal{G}_{\ell^\infty}^{sc}(\Q_\ell)$ and the isogeny $p^{sc}:\mathcal{G}_{\ell^\infty }^{sc}\rightarrow \mathcal{G}_{\ell^\infty }$ extends uniquely to a  morphism  $p^{sc}:\frak{G}_{\ell^\infty }^{sc}\rightarrow \frak{G}_{\ell^\infty }$ of group schemes over $\Z_\ell$. Furthermore,  $p^{sc}:\mathcal{G}_{\ell^\infty }^{sc}\rightarrow \mathcal{G}_{\ell^\infty }$ induces a morphism  (Fact \ref{Fact:Frattini})
$ \frak{G}_{\ell^\infty }^{sc}(\Z_\ell)\rightarrow \frak{G}_{\ell^\infty }(\Z_\ell)^+$ with the following properties.
\begin{itemize}
\item (i) The diagram 
$$\xymatrix{\mathcal{G}_{\ell^\infty }^{sc}(\Q_\ell)\ar[r]^{p^{sc}}& \mathcal{G}_{\ell^\infty }(\Q_\ell)\\
 \frak{G}_{\ell^\infty }^{sc}(\Z_\ell)\ar[r]\ar@{_{(}->}[u]&\frak{G}_{\ell^\infty }(\Z_\ell)^+\ar@{_{(}->}[u]}$$
is cartesian. Indeed, since  $p^{sc}:\mathcal{G}_{\ell^\infty }^{sc}\rightarrow \mathcal{G}_{\ell^\infty }$ is open (even a local isomorphism), $\frak{G}_{\ell^\infty }(\Z_\ell)^+\cap p^{sc}(\mathcal{G}_{\ell^\infty}^{sc}(\Q_\ell))\subset \frak{G}_{\ell^\infty }(\Z_\ell)$ is open  hence closed, hence compact. So $(p^{sc})^{-1}(\frak{G}_{\ell^\infty }(\Z_\ell)^+)$ is compact in $\mathcal{G}_{\ell^\infty }^{sc}(\Q_\ell)$ as an extension of the compact group $\frak{G}_{\ell^\infty }(\Z_\ell)^+\cap p^{sc}(\mathcal{G}_{\ell^\infty}^{sc}(\Q_\ell))$ by a finite group. Also,  $(p^{sc})^{-1}(\frak{G}_{\ell^\infty }(\Z_\ell)^+)$ contains $\frak{G}_{\ell^\infty }^{sc}(\Z_\ell)$. Thus, by maximality of $\frak{G}_{\ell^\infty }^{sc}(\Z_\ell)$, we have $(p^{sc})^{-1}(\frak{G}_{\ell^\infty }(\Z_\ell)^+) =\frak{G}^{sc}_{\ell^\infty }(\Z_\ell)$.\\
\item (ii) The homomorphism $ \frak{G}_{\ell^\infty }^{sc}(\Z_\ell)\rightarrow \frak{G}_{\ell^\infty }(\Z_\ell)^+$ is surjective. Indeed, let $g\in \frak{G}_{\ell^\infty }(\Z_\ell)^+$. Since the image of $p^{sc}:\mathcal{G}_{\ell^\infty }^{sc}(\Q_\ell)\rightarrow \mathcal{G}_{\ell^\infty }(\Q_\ell)$ is normal and its cokernel is of  exponent bounded by a constant depending only on the rank of $\mathcal{G}_{\ell^\infty }$, $g$ lies in the image of $p^{sc}:\mathcal{G}_{\ell^\infty }^{sc}(\Q_\ell)\rightarrow \mathcal{G}_{\ell^\infty }(\Q_\ell)$  for $\ell\gg 0$ (compared with the rank of  $\mathcal{G}_{\ell^\infty }$), that is, by (i), in the image of $ (p^{sc})^{-1}(\frak{G}_{\ell^\infty }(\Z_\ell)^+)= \frak{G}_{\ell^\infty }^{sc}(\Z_\ell)$. 
\end{itemize}
(\ref{Cor:SSConsequence3}) thus follows from (\ref{Cor:SSConsequence1}.2).\\  
\indent  The implication (\ref{Cor:SSConsequence3}) $\Rightarrow$ (\ref{Cor:SSConsequence2}.4) follows from Lemma \ref{Lem:Compare} and (\ref{Cor:SSConsequence1}.2). $\square$\\

\noindent This concludes the proof of Theorem \ref{Th:GSSMax}.\\

 \begin{center}\textbf{\fontfamily{lmr}\selectfont{PART III: SEMISIMPLICITY}}\end{center}
  
  \section{Lie-theoretic proof of Theorem \ref{Th:GSS}}\label{Proof1}
  The results explained here are entirely due to the second author, Chun-Yin Hui. They led to the first complete proof of Theorem \ref{Th:GSS}.   
  \subsection{Semisimple models and good lattices} Let $\mathcal{G}_{\ell^\infty}$ be a connected semisimple group over $\Q_\ell$ of dimension $\delta$ and rank $s$. Let $V_{\ell^\infty}$ be a faithful, $r$-dimensional $\Q_\ell$-representation of $\mathcal{G}_{\ell^\infty}$. Fix a  lattice $H_{\ell^\infty}\hookrightarrow V_{\ell^\infty}$; this defines a model $\SGL_{H_{\ell^\infty}}$ of $\SGL_{V_{\ell^\infty}}$ over $\Z_\ell$. Let $\frak{G}_{\ell^\infty}$ denote the  Zariski closure of $\mathcal{G}_{\ell^\infty} $ inside $\SGL_{H_{\ell^\infty}}$ (endowed with the reduced subscheme structure). Then   $\frak{G}_{ \ell^\infty}$ is a flat model of $\mathcal{G}_{\ell^\infty}$ over $\Z_\ell$. Under mild assumptions, we give a criterion in terms of  tensor-invariants data to ensure that $\frak{G}_{ \ell^\infty}$ is a semisimple group scheme over $\Z_\ell$. Write $\mathcal{G}_{\ell}:=\frak{G}_{ \ell^\infty,\F_\ell}$ and $\mathcal{G}_{\ell}^\circ$ for its identity component.\\
  
  \noindent Let $\hbox{\rm Rep}^f_{\Z_\ell}(\frak{G}_{\ell^\infty})$ denote the category of  finitely generated free $\Z_\ell$-modules $M$ together with a morphism of $\Z_\ell$-group schemes $\frak{G}_{\ell^\infty}\rightarrow \SGL_M$. Define 
  $$\begin{tabular}[t]{lcll}
  $\Delta_{H_{\ell^\infty}}: $&$\hbox{\rm Rep}^f_{\Z_\ell}(\frak{G}_{\ell^\infty})$&$\rightarrow$&$ \Z_{\geq 0}$\\
  &$M$&$\mapsto$&$\Sdim_{\F_\ell}(M_{\F_\ell}^{\mathcal{G}_\ell^\circ})-\Sdim_{\Q_\ell}(M_{\Q_\ell}^{\mathcal{G}_{\ell^\infty}})$.
  \end{tabular}$$
  
  \noindent Let  $\mathcal{T}_{\ell^\infty}\subset \mathcal{G}_{\ell^\infty}$ be a maximal torus. We will say that $\mathcal{T}_{\ell^\infty}$ admits a nice model with respect to  $H_{\ell^\infty}$  if  $\mathcal{T}_{\ell^\infty}$ splits over a finite
  extension  $E_\ell$ 
  of $\Q_\ell$  and  if   the closed embedding $\mathcal{T}_{\ell^\infty,E_\ell}\simeq \G_{m,E_\ell}^s\hookrightarrow \SGL_{V_{\ell^\infty,E_\ell}}$ extends to a closed embedding $\G_{m,\mathcal{O}_\ell}^s\hookrightarrow\SGL_{H_{\ell^\infty,\mathcal{O}_\ell}}$, where $\mathcal{O}_\ell$ denotes the ring of integers of $E_\ell$. 
 
  \begin{theoreme}\label{Th:Lie} Assume  $\mathcal{G}_{\ell^\infty}$ contains a maximal torus which  admits a nice model with respect to  $H_{\ell^\infty}$. Then,
    \begin{itemize}[leftmargin=* ,parsep=0cm,itemsep=0cm,topsep=0cm]
  \item (\ref{Th:Lie}.1) $\frak{G}_{ \ell^\infty}$ is smooth over $\Z_\ell$;
  \item (\ref{Th:Lie}.2)  The quotient $\mathcal{G}_\ell^{rd}$ of $\mathcal{G}_\ell^\circ$ by its unipotent radical is a reductive group of rank $s$ (and the root system of $\mathcal{G}_{\ell,\overline{\F}_\ell}^{rd}$ is a subsystem of the root system of  $\mathcal{G}_{\ell^\infty,\overline{\Q}_\ell}$);
  \item (\ref{Th:Lie}.3)  For $\ell\gg 0 $ (depending only on $r$) the following holds. Let $\frak{g}_{\ell^\infty}\hookrightarrow H_{\ell^\infty}\otimes H_{\ell^\infty}^\vee$ denote the Lie algebra of $\frak{G}_{\ell^\infty}$. Then  $\frak{G}_{\ell^\infty}$ is semisimple over $\Z_\ell$ if and only if $\Delta_{ H_{\ell^\infty}}(M)\leq 0$ for $M=\Lambda^n\frak{g}_{\ell^\infty}^\vee$, $n=1,\dots, \delta$.  
  \end{itemize}
  \end{theoreme}
  \noindent\textit{Proof.} As $\spec(\mathcal{O}_\ell)\rightarrow\spec(\Z_\ell)$ is flat, $ \frak{G}_{ \ell^\infty, \mathcal{O}_\ell}\subset \SGL_{H_{\ell^\infty, \mathcal{O}_\ell}}$   coincides with the Zariski closure of $\mathcal{G}_{\ell^\infty, E_\ell}$ (endowed with its reduced subscheme structure) and $\mathcal{G}_{\ell^\infty , E_\ell}=\frak{G}_{\ell^\infty, E_\ell}$ \cite[1.2.6]{BT84}. So to perform the proof, we may base-change to $\mathcal{O}_\ell$, $E_\ell$. For simplicity, we assume $\mathcal{O}_\ell=\Z_\ell$ and $E_\ell=\Q_\ell$ below.\\
 \indent Fix a Borel subgroup $\mathcal{T}_{\ell^\infty }\subset \mathcal{B}_{\ell^\infty}\subset \mathcal{G}_{\ell^\infty}$; write $\Phi:=\Phi(\mathcal{G}_{\ell^\infty},\mathcal{T}_{\ell^\infty })$ for the root system and let $\Phi^+\subset \Phi$ denote the set of positive roots defined by $\mathcal{B}_{\ell^\infty}$. For $\alpha\in \Phi$, let $\frak{g}_{\ell^\infty,\alpha}\subset \frak{g}_{\ell^\infty,\Q_\ell} =\hbox{\rm Lie}(\mathcal{G}_{\ell^\infty} )\subset \frak{g}\frak{l}(V_{\ell^\infty} )$ and $\mathcal{U}_{\ell^\infty,\alpha}\subset \mathcal{G}_{\ell^\infty }$ denote the corresponding root space and group respectively. Let $\frak{T}_{\ell^\infty}\simeq \G_{m,\Z_\ell}^s, \frak{U}_{\ell^\infty,\alpha}\subset \frak{G}_{\ell^\infty }$ denote the Zariski closure of   $\mathcal{T}_{\ell^\infty}, \mathcal{U}_{\ell^\infty,\alpha}$; by construction $\frak{T}_{\ell^\infty}, \frak{U}_{\ell^\infty,\alpha}$ are flat over $\Z_\ell$. \\ 
 \indent For $\alpha\in \Phi $,  $\frak{g}_{ \ell^\infty,\alpha}\cap \frak{g}\frak{l}(H_{\ell^\infty})$ is a free $\Z_\ell$-module of rank $1$. Let $N_{\ell^\infty,\alpha}\in \frak{g}_{ \ell^\infty,\alpha}\cap \frak{g}\frak{l}(H_{\ell^\infty})$ be a $\Z_\ell$-basis (in particular $N_{\ell^\infty,\alpha}\otimes\F_\ell\not=0$). Then for $\ell\geq r$, the closed embedding $x_{\ell^\infty,\alpha}: \G_{a,\Z_\ell}\rightarrow \frak{G}_{\ell^\infty}$, $ t\mapsto \exp(tN_{\ell^\infty,\alpha})$ induces an isomorphism of $\Z_\ell$-group schemes onto a closed subgroup scheme of $\frak{G}_{ \ell^\infty }$ which coincides with
  $\frak{U}_{\ell^\infty,\alpha}\subset \frak{G}_{ \ell^\infty}$.\\
  \indent By \cite[2.2.3 (iii)]{BT84} the $\frak{T}_{\ell^\infty}$-equivariant morphism induced by multiplication
$$ \prod_{\alpha\in \Phi^+}\frak{U}_{\ell^\infty,\alpha}\times\frak{T}_{\ell^\infty}\times  \prod_{\alpha \in -\Phi^+}  \frak{U}_{\ell^\infty,\alpha}\rightarrow \frak{G}_{\ell^\infty}$$
induces an isomorphism onto a dense open $\Z_\ell$-subscheme of $ \frak{G}_{ \ell^\infty}$. In particular $\frak{G}_{\ell^\infty }$  is smooth over $\Z_\ell$, $\frak{G}_{\ell^\infty}^\circ\subset \frak{G}_{\ell^\infty}$ is an open subgroup scheme, smooth, affine over $\Z_\ell$ \cite[2.2.5]{BT84} and   $\mathcal{G}_{\ell}^\circ$ contains a split torus of rank $s$. As the reductive rank is lower semicontinous (\textit{e.g.} \cite[Thm. 10.4.2]{Gille}), $\mathcal{G}_{\ell}^\circ$ has reductive rank $s$. Furthermore, reduction modulo-$\ell$ induces a canonical isomorphism of $\G_m^s$-modules $\frak{g}_{ \ell^\infty } \otimes\F_\ell\tilde{\rightarrow} \hbox{\rm Lie}(\mathcal{G}_{ \ell})$. The latter implies that the root system of $\mathcal{G}^{rd}_{\ell , \overline{\F}_\ell}$ is a subsystem of the root system of  $\mathcal{G}_{\ell^\infty ,\overline{\Q}_\ell}$ (since $\hbox{\rm Lie}(\mathcal{G}^{rd}_{\ell})$ is a quotient of $\hbox{\rm Lie}(\mathcal{G}_{ \ell})$).\\
\indent We now turn to the proof of (\ref{Th:Lie}.3). Let $\cG^u_{ \ell}$ denote the unipotent radical of $\cG_{ \ell}$ and write $\frak{g}_\ell:=\hbox{\rm Lie}(\mathcal{G}_{ \ell})$,  $\frak{g}_\ell^u:=\hbox{\rm Lie}(\mathcal{G}_{ \ell}^u)$, $\frak{g}_\ell^{rd}:=\hbox{\rm Lie}(\mathcal{G}_{ \ell}^{rd})$. This gives rise to a decomposition of the adjoint representation  
 
$$0\to \frak{g}_\ell^u \to \frak{g}_\ell\to \frak{g}_\ell^{rd}\to 0,$$
and dualizing, 
$$
0\to (\frak{g}_\ell^{rd})^\vee\to \frak{g}_\ell^\vee\to (\frak{g}_\ell^{u})^\vee\to 0.
$$
As  $\cG^u_{\ell}$ acts trivially on $\frak{g}_\ell^{rd}$ (observe that for  $g\in \cG^u_{\ell}$ the conjugation automorphism $c_g$ on $\cG_{ \ell}^\circ$ descends to the 
identity map on $\cG^{rd}_{ \ell}=\cG_{ \ell}^\circ/\cG^u_{\ell}$), the action of $\mathcal{G}_\ell^\circ$ on $\frak{g}_\ell^{rd}$  factors through the adjoint representation of $\mathcal{G}_\ell^{rd}$.\\

\noindent Suppose now the condition on
$\Delta_{H_{\ell^\infty}}$ is satisfied.\\

\noindent\textbf{Claim 1:} \textit{For $\ell\gg 0 $ (depending only on $r$), $\mathcal{G}_\ell^{rd}$ is semisimple.}\\

\noindent \textit{Proof of Claim 1.} Compute the dimension of the Lie algebra $\frak{z}_\ell$ of the center of $\mathcal{G}_\ell^{rd}$ $$\begin{tabular}[t]{lll}
$\Sdim_{\F_\ell}(\frak{z}_\ell)$&$\leq \Sdim_{\F_\ell}((\frak{g}_{\ell }^{rd})^{\mathcal{G}_\ell^{rd}})\stackrel{(1)}{=}\Sdim_{\F_\ell}((\frak{g}_{\ell }^{rd\vee})^{\mathcal{G}_\ell^{rd}})$&$=\Sdim_{\F_\ell}((\frak{g}_{\ell }^{rd\vee})^{\mathcal{G}^\circ_\ell})$\\
&&$\stackrel{(2)}{\leq} \Sdim_{\F_\ell}((\frak{g}_{\ell }^\vee)^{\mathcal{G}^\circ_\ell})$\\
&&$\stackrel{(3)}{\leq} \Sdim_{\Q_\ell}((\frak{g}_{\ell^\infty,\Q_\ell}^\vee)^{\mathcal{G}_{\ell^\infty}})\stackrel{(1)}{=}\Sdim_{\Q_\ell}(\frak{g}_{\ell^\infty,\Q_\ell}^{\mathcal{G}_{\ell^\infty}})\stackrel{(4)}{=}0,$
\end{tabular}$$
where (1) is by the semisimplicity  \cite[Prop. 3.2]{Jantzen} (see also \cite[Cor. 4.3]{Springer}) and self-duality of the adjoint representation of the reductive group $\mathcal{G}_\ell^{rd}$ (resp. $\mathcal{G}_{\ell^\infty}$) for $\ell\gg 0$ compared with the rank (resp. for all $\ell$), (2) is because $\frak{g}_{\ell }^{rd\vee}$ is a submodule of $\frak{g}_{\ell }^\vee$, (3) is the assumption $\Delta_{H_{\ell^\infty}}(\frak{g}_{\ell^\infty}^\vee)\leq 0$ and  (4) is because $\mathcal{G}_{\ell^\infty}$ is semisimple. $\square$\\

\noindent\textbf{Claim 2:}  \textit{For every integer $n$ and for $\ell\gg 0$ (depending only on $n$) the following holds. Consider a pair of rank $n$ connected semisimple groups $\mathcal{G}$ over $\overline{\F}_\ell$ and $\mathcal{G}'$ over $\overline{\Q}_\ell$. Assume $\dim(\mathcal{G})<\dim(\mathcal{G}')$. Then there exists $0\leq m\leq \Sdim(\mathcal{G})$ such that $\dim ((\Lambda^m \mathfrak{g})^{\mathcal{G}})> \dim ((\Lambda^m \mathfrak{g}')^{\mathcal{G}'})$. Here $\frak{g}$ and $\frak{g}'$ denote the Lie algebra of $\mathcal{G}$ and $\mathcal{G}'$ respectively. }\\ 

\noindent\textit{Proof of Claim 2.} The assertion will follow from the explicit computation of the invariant dimensions of the exterior algebra $\Lambda^*\mathfrak{g}$ for a rank $n$ connected semisimple algebraic group $\mathcal{G}$ over an algebraically closed field. Assume first $\mathcal{G}$ is almost simple. Over $\C$ (and thus over any algebraically closed field of characteristic zero), these  are given by the coefficients (corresponding to the exterior powers) of the Poincar\'e polynomial $P_{\mathcal{G}}(T)$ of the cohomology of the Lie group $\mathcal{G}$ (\textit{e.g.}  \cite[$\mathsection0$]{Ba01}):
$$ P_{\mathcal{G}}(T)=\prod_{j=1}^n(1+T^{2d_j+1}),$$
\noindent where $d_1,...,d_n$ are the exponents of the Weyl group (for the explicit values of $d_1,\dots ,d_n$ for each simple type, see\cite[Prop. 10.2.5]{Ca72}). These results also hold for a connected almost simple algebraic group $\mathcal{G}$ over $\bar\F_\ell$ when $\ell\gg 0$ compared with the rank of $\mathcal{G}$ as can be shown from the semisimplicity of the representations  \cite[Prop. 3.2]{Jantzen} (see also \cite[Cor. 4.3]{Springer} and \cite{SerreSS}) and the classification of irreducible representations of the simply connected cover $\mathcal{G}^{sc}$ \cite[$\mathsection12$ Thm. 41]{St67}.  For an arbitrary connected semisimple group $\mathcal{G}$ over an algebraically closed field,  the (graded) invariant dimensions of $\Lambda^*\mathfrak{g}$ are given by the coefficients of 
the product of Poincar\'e polynomials $$P_{\mathcal{G}}(T)=\prod_{i=1}^t P_{\mathcal{G}_i}(T),$$ where $\mathcal{G}_1,...,\mathcal{G}_t$ are the  almost simple factors of $\mathcal{G}$. This follows easily from the algebra isomorphism  $\Lambda^*(V\oplus W)=\Lambda^*(V)\otimes\Lambda^*(W)$.\\
\indent Now, we apply the above to $\mathcal{G}$ and $\mathcal{G}'$ as in Claim 2. On the one hand, we have
$$P_{\mathcal{G}}(1)=2^n=P_{\mathcal{G}'}(1)$$ 
while, on the other hand 
$$\Sdeg(P_{\mathcal{G}})=\Sdim(\mathcal{G})<\Sdim(\mathcal{G}')=\Sdeg(P_{\mathcal{G}'}).$$
This shows that the sum of coefficients in degrees $\leq \Sdim(\mathcal{G})$ of $P_{\mathcal{G}}$ is strictly larger than  the sum of coefficients in degrees $\leq \Sdim(\mathcal{G})$ of $P_{\mathcal{G}'}$. In particular, there exists $m\leq  \Sdim(\mathcal{G})$ such that the coefficient in degree $m$ of $P_{\mathcal{G}}$ is strictly larger than the  coefficient in degree $m$ of $P_{\mathcal{G}'}$.
$\square$\\

    \noindent We can now conclude the proof of (\ref{Th:Lie}.3). Let us first show that the condition on $\Delta_{H_{\ell^\infty}}$ is sufficient.  As $\frak{G}_{\ell^\infty}$ is flat over $\Z_\ell$, it is enough to show that (i)  $\mathcal{G}_{\ell}$ is connected and (ii) $\mathcal{G}_\ell^{rd} $  and $\mathcal{G}_{\ell^\infty}$ have the same dimension. Assertion (i) follows from (ii) and the fact that $\mathcal{G}_{\ell^\infty}$ is connected \cite[Prop. 3.1.3]{Conrad}. If (ii) does not hold, then   there exists  $ 0\leq m\leq \delta$ such that   $\dim_{\F_\ell} ((\Lambda^m \mathfrak{g}_\ell^{rd})^{\mathcal{G}_\ell^{rd}})> \dim_{\Q_\ell} ((\Lambda^m \mathfrak{g}_{\ell^\infty,\Q_\ell})^{\mathcal{G}_{\ell^\infty}})$ by Claim 2 applied to $\mathcal{G}=\mathcal{G}_{\ell,\overline{\F}_\ell}^{rd}$ (which is semisimple by Claim 1) and 
$\mathcal{G}'=\mathcal{G}_{\ell^\infty,\overline{\Q}_\ell}$ and the rank assertion in (\ref{Th:Lie}.2). Since $\Delta_{H_{\ell^\infty}}(\Lambda^m \mathfrak{g}_\ell^\vee)\leq 0$, this contradicts the following inequalities (for $\ell\gg 0$):
$$\begin{tabular}[t]{ll} $\Sdim_{\F_\ell}((\Lambda^m\frak{g}_{\ell }^{rd})^{\mathcal{G}_\ell^{rd}})$&$=\Sdim_{\F_\ell}((\Lambda^m\frak{g}_{\ell }^{rd\vee})^{\mathcal{G}_\ell^{rd}}) $\\
&$=\Sdim_{\F_\ell}((\Lambda^m\frak{g}_{\ell }^{rd\vee})^{\mathcal{G}^\circ_\ell})$\\
&$ \leq  \Sdim_{\F_\ell}((\Lambda^m\frak{g}_{\ell }^\vee)^{\mathcal{G}^\circ_\ell})$\\
&$\leq \Sdim_{\Q_\ell}((\Lambda^m\frak{g}_{\ell^\infty,\Q_\ell}^\vee)^{\mathcal{G}_{\ell^\infty}}) =\Sdim_{\Q_\ell}((\Lambda^m\frak{g}_{\ell^\infty,\Q_\ell})^{\mathcal{G}_{\ell^\infty}}).$\\
\end{tabular}$$
The argument also shows that the  condition on $\Delta_{H_{\ell^\infty}}$ is necessary since the root system  (hence, the Poincar\'e polynomial) of a semisimple group scheme is locally constant. $\square$

\subsection{Application to the Proof of Theorem \ref{Th:GSS}} We retain the notation of Section \ref{Summary} and Part II. From Corollary \ref{Cor:SSConsequence2}, it is enough to prove that $\frak{G}_{\ell^\infty}$ is a semisimple group scheme over $\Z_\ell$. This can be checked by applying the criterion of Theorem \ref{Th:Lie}.

\subsubsection{}\label{Verif} The fact that the assumptions of Theorem \ref{Th:Lie} are satisfied for $\ell\gg 0$ follows from the first paragraph in the proof of \cite[Prop. 1.3]{LarsenPinkAV}. Indeed, one can always find a $\Gamma$-regular element \footnote{\label{GReg}Let $Q$ be a field of characteristic $0$, $V$ a finite-dimensional $Q$-vector space and $G\subset \SGL_V$   a reductive subgroup. Then a regular semisimple element $g\in G(Q)$ is said to be $\Gamma$-regular for $G\subset \SGL_V$ if  every automorphism of $T_g\times_Q\overline{Q}$ which fixes $g$ and preserves the formal character of $T_g\subset \SGL_V$ is trivial, and if the only $\SGL_V(\overline{Q})$-conjugate of $T_g\times_Q\overline{Q}$ containing $g$ is $T_g\times_Q\overline{Q}$. Here $T_g\subset G$ denotes the (necessarily unique) maximal torus containing $g$. We refer to \cite[\S 1]{LarsenPinkAV} and \cite[\S 4]{LP} for details.} $t$ for $ \mathcal{G}_{0\ell^\infty} \subset \SGL_{V_{\ell^\infty}}$ with the property that the characteristic polynomial of $t$ acting on $V_{\ell^\infty}$ coincides with the characteristic polynomial $P_{x_0}$ of $\rho_{\ell^\infty}(F_{x_0})$ for some $x_0\in X_0$. More precisely, let $V_{\ell^\infty}^{ss}$ denote the $\pi_1(X_0,x)$-semisimplification of $V_{\ell^{\infty}}$ and $\mathcal{G}_{0\ell^\infty}^{rd}$ the Zariski closure of the image of $\pi_1(X_0,x)$ acting on  $V_{\ell^\infty}^{ss}$. Note that $\mathcal{G}_{0\ell^\infty}^{rd}$ identifies with the quotient of $\mathcal{G}_{0\ell^\infty}$ by its unipotent radical. In particular (\ref{Summary}.2.2), $\mathcal{G}_{0\ell^\infty}^{rd}$ is connected reductive. By \cite[Prop. 7.2]{LP}, there exists (a density $1$ set of) $x_0\in \pi_1(X_0,x)$ such that the image $t_{x_0}$ of $F_{x_0}$ by $\pi_1(X_0,x)\rightarrow\SGL(V_{\ell^\infty}^{ss})$ is $\Gamma$-regular for $ \mathcal{G}_{0\ell^\infty}^{rd}\subset \SGL_{V_{\ell^\infty}^{ss}}$. Let $\mathcal{T}_{0\ell^\infty}^{rd}\subset \mathcal{G}_{0\ell^\infty}^{rd}$ denote the corresponding maximal torus. Since the kernel of $ \mathcal{G}_{0\ell^\infty} \twoheadrightarrow \mathcal{G}_{0\ell^\infty}^{rd}$ is the unipotent radical of $\mathcal{G}_{0\ell^\infty}$, there exists a  maximal torus $\mathcal{T}_{0\ell^\infty}\subset \mathcal{G}_{0\ell^\infty} $ lifting $\mathcal{T}_{0\ell^\infty}^{rd}$   and mapping isomorphically onto $\mathcal{T}_{0\ell^\infty}^{rd}$. Then the unique element $t\in \mathcal{T}_{0\ell^\infty}(\Q_\ell)$ lifting $t_{x_0}$ has the desired property.\\

\indent Let $\mathcal{T}_{0\ell^\infty} \subset \mathcal{G}_{0\ell^\infty} $ denote the unique (necessarily maximal) torus containing $t$. Let $\mathcal{T}_{\ell^\infty} \subset \mathcal{G}_{\ell^\infty} $ denote the maximal torus of $\mathcal{G}_{\ell^\infty} $ contained in $\mathcal{T}_{0\ell^\infty}$. By definition of $\Gamma$-regularity (see Footnote \ref{GReg}), the splitting fields $E_\ell/\Q_\ell$ of $\mathcal{T}_{0\ell^\infty} $ and  $P_{x_0}$ over $\Q_\ell$ coincide. In particular, for $\ell\gg 0$ (not dividing the discriminant of the product of the monic irreducible factors of $P_{x_0}$) the eigenspace decomposition $V_{\ell^\infty,E_\ell}=\oplus_\lambda V_{\ell^\infty,E_\ell}(\lambda)$ of $t$ coincides with the one of $\mathcal{T}_{0\ell^\infty,E_\ell}$ and induces a decomposition $H_{\ell^\infty,\mathcal{O}_\ell}=\oplus_\lambda (H_{\ell^\infty,\mathcal{O}_\ell}\cap V_{\ell^\infty, E_\ell}(\lambda))$. This ensures that the closed embedding $\mathcal{T}_{\ell^\infty,E_\ell}\simeq \G_{m,E_\ell}^s\hookrightarrow \SGL_{V_{\ell^\infty,E_\ell}}$ extends to a closed embedding $\G_{m,\mathcal{O}_\ell}^s\hookrightarrow\SGL_{H_{\ell^\infty,\mathcal{O}_\ell}}$.

\subsubsection{}\label{Last} Let $\delta$ denote the dimension of $\mathcal{G}_{\ell^\infty}$.  It only remains to show that, for $\ell\gg 0$,  $\Delta_{H_{\ell^\infty}}(M)\leq 0$ for $M=\Lambda^n\frak{g}_{\ell^\infty}^\vee$, $n=1,\dots, \delta$. Note that, from Theorem   \ref{Cor:SSConsequence1}, $\mathcal{G}_\ell$ is connected. Also, as $\Pi_{\ell^\infty}$ is normal in $\Pi_{0\ell^\infty}$, $\frak{g}_{\ell^\infty}$ is a $\Pi_{0\ell^\infty}$-module. \\

\noindent From (\ref{Summary}.3.1) applied to the $\Pi_{0\ell^\infty}$-module  quotients
$  H_{\ell^\infty}^{\otimes n}\otimes (H_{\ell^\infty}^\vee)^{\otimes n}\twoheadrightarrow \Lambda^n\frak{g}_{\ell^\infty}^\vee$ it is enough to show that 
for  $M=\Lambda^n\frak{g}_{\ell^\infty}^\vee$, $n=1,\dots, \delta$ 
we have $$(\ref{Last}.1)\;\; \Sdim(M_{\F_\ell}^{\mathcal{G}_\ell})\leq \Sdim(M_{\F_\ell}^{\Pi_\ell}) 
\;\; \hbox{\rm and}\;\;(\ref{Last}.2)\;\; \Sdim(M_{\Q_\ell}^{\mathcal{G}_{\ell^\infty}})=\Sdim(M_{\Q_\ell}^{\Pi_{\ell^\infty}}).$$
\noindent  (\ref{Last}.1) always holds since $\Pi_\ell\subset \mathcal{G}_\ell(\F_\ell)$. As for (\ref{Last}.2), for every  $ v\in M_{\Q_\ell}^{\Pi_{\ell^\infty}}$, we have $\Pi_{\ell^\infty}\subset \hbox{\rm Stab}_{\hbox{\rm \tiny GL}_{V_{\ell^\infty}}}(v)=:\mathcal{S}_v$. But as $\mathcal{S}_v\subset \SGL_{V_{\ell^\infty}}$ is a closed algebraic subgroup, this implies $\mathcal{G}_{\ell^\infty}\subset \mathcal{S}_v$. \\

\section{Cohomological proof of Theorem \ref{Th:GSS}}\label{Proof2}
\subsection{Preliminary reductions}
\begin{lemme}\label{Lem:Red}
 Let $\iota:A_\ell\hookrightarrow H_\ell$ be a $\Pi_\ell$-submodule of $\F_\ell$-dimension $a$. Consider the following conditions.
$$\begin{tabular}[t]{ll}
(\ref{Lem:Red}.1)& $\iota: A_\ell\hookrightarrow H_\ell$ splits as a morphism of $\Pi_\ell$-modules; \\
 (\ref{Lem:Red}.2)& $\Lambda^a\iota: \Lambda^aA_\ell\hookrightarrow \Lambda^aH_\ell$ splits as a morphism of $\Pi_\ell$-modules;\\
 (\ref{Lem:Red}.3)& $ (\Lambda^aH_\ell)^{\Pi_\ell}\hookrightarrow \Lambda^aH_\ell$ splits as a morphism of $\Pi_\ell$-modules.
 \end{tabular}$$
 Then, for $\ell\gg 0$, (\ref{Lem:Red}.3) $\Rightarrow$  (\ref{Lem:Red}.2) $\Rightarrow$  (\ref{Lem:Red}.1).
 \end{lemme}
\noindent\textit{Proof.} (\ref{Lem:Red}.3) $\Rightarrow$ (\ref{Lem:Red}.2): As $\Pi_\ell=\Pi_\ell^+$ (\ref{Summary}.2.1), $\Pi_\ell$ acts trivially on $ \Lambda^aA_\ell$ that is $\Lambda^aA_\ell\hookrightarrow (\Lambda^aH_\ell)^{\Pi_\ell}$ and this trivially splits as a morphism of $\Pi_\ell$-modules. \\
\noindent (\ref{Lem:Red}.2) $\Rightarrow$ (\ref{Lem:Red}.1): Fix a $\Pi_\ell$-equivariant splitting $s:\Lambda^aH_\ell\rightarrow \Lambda^aA_\ell$ of $\Lambda^a\iota: \Lambda^aA_\ell\hookrightarrow \Lambda^aH_\ell$. Then one can explicitly construct  a $\Pi_\ell$-equivariant splitting for $\iota:  A_\ell\hookrightarrow H_\ell$ as follows
$$\xymatrix{H_\ell\ar[rr]^{v\rightarrow v\wedge-}\ar@{.>}[dr]^{\sigma_s}&&\SHom(\Lambda^{a-1}H_\ell,\Lambda^aH_\ell)\ar[r]^{-\circ\wedge^{a-1}\iota}&\SHom(\Lambda^{a-1}A_\ell,\Lambda^aH_\ell)\ar[d]^{s\circ-}\\
&A_\ell\ar[rr]_{\simeq}^{v\rightarrow v\wedge-}&&\SHom(\Lambda^{a-1}A_\ell,\Lambda^aA_\ell)\; .\square}$$

\noindent From Lemma \ref{Lem:Red}, it is enough to show (\ref{Lem:Red}.3) for $H_\ell=\SH^w(Y_x,\F_\ell)$. For $\ell\gg 0$,    $\Lambda^aH_\ell$ is a direct factor of   $H^{wa}(Y_x^{[a]},\F_\ell)$. So replacing $Y_0\rightarrow X_0$ with $Y_0^{[a]}\rightarrow X_0$, it is enough to show that $H_\ell^{\Pi_{\ell }}\hookrightarrow H_\ell$ splits as a morphism of $\Pi_\ell$-modules. That is, writing $A_\ell:= H_\ell^{\Pi_{\ell }}$ and $B_\ell:=H_\ell/ H_\ell^{\Pi_{\ell }}$, 
$$\begin{tabular}[t]{ll}
(\ref{Lem:Red}.4)& The class $[H_\ell]$ of $0\rightarrow A_\ell \rightarrow H_\ell\rightarrow  B_\ell \rightarrow 0$
is trivial in $H^1(\Pi_{\ell}, A_\ell\otimes B_\ell^\vee)$ (or, equivalently, in\\
& $H^1(\Pi_{\ell^\infty}, A_\ell\otimes B_\ell^\vee)$).
\end{tabular}$$

 \subsection{Proof of (\ref{Lem:Red}.4)}\label{Pr} Write $B_{\ell^\infty}:=H_{\ell^\infty}/H_{\ell^\infty}^{\Pi_{\ell^\infty}}$, $A_{\ell^\infty}:=H_{\ell^\infty}^{\Pi_{\ell^\infty}}$ . \\

\indent Consider the following exact commutative snake diagram  $$\xymatrix{&0\ar[d]&0\ar[d]&0\ar[d]&\\
0\ar[r]&A_{\ell^\infty} \ar[r]^\ell\ar@{.>}[d]&A_{\ell^\infty} \ar[r]\ar[d]&A_{\ell } \ar[r]\ar[d]&C'\ar[r]\ar@{.>}[d]&0\\
0\ar[r]&H_{\ell^\infty} \ar[r]^\ell\ar@{.>}[d]&H_{\ell^\infty}\ar[r]\ar[d]&H_{\ell} \ar[r]\ar[d]&C\ar[r]\ar@{.>}[d]&0\\
0\ar[r]&K''\ar[r]^\iota \ar@{.>}[rrruu]^\delta&B_{\ell^\infty} \ar[r]\ar[d]&B_{\ell} \ar[r]\ar[d]&C''\ar[r]\ar[d]&0\\
&&0&0&0&}$$
From Fact \ref{Gabber} $C=0$ hence $C''=0$ while from (\ref{Summary}.3), $C'=0$. This shows  $\delta=0$ hence that  there is an isomorphism $K''\simeq B_{\ell^\infty}$ with respect to which $\iota$ identifies with multiplication-by-$\ell$. In particular,\\

\begin{itemize}[leftmargin=* ,parsep=0cm,itemsep=0cm,topsep=0cm]
\item (\ref{Pr}.1) the sequence $ 0\rightarrow B_{\ell^\infty}\stackrel{\ell}{\rightarrow} B_{\ell^\infty}\rightarrow B_\ell\rightarrow 0$
is exact.\\
\end{itemize}
\indent   Let  $[H_{\ell^\infty}]$ and $[H_{\ell^\infty \Q_\ell}]$ denote the class of the extensions
 $$0\rightarrow  A_{\ell^\infty}\rightarrow H_{\ell^\infty}\rightarrow   B_{\ell^\infty}\rightarrow 0$$
 $$0\rightarrow  A_{\ell^\infty \Q_\ell}\rightarrow H_{\ell^\infty \Q_\ell}\rightarrow   B_{\ell^\infty \Q_\ell}\rightarrow 0$$
 in  $H^1(\Pi_{\ell^\infty}, A_{\ell^\infty}\otimes B_{\ell^\infty}^{\vee})$, $H^1(\Pi_{\ell^\infty}, A_{\ell^\infty \Q_\ell}\otimes B_{\ell^\infty \Q_\ell}^\vee )$ respectively. \\
 \noindent Then $[H_{\ell^\infty}]$ maps to $[H_\ell]$ \textit{via} $$H^1(\Pi_{\ell^\infty}, A_{\ell^\infty}\otimes B_{\ell^\infty}^{\vee})\rightarrow H^1(\Pi_{\ell^\infty}, A_\ell\otimes B_\ell^\vee)$$ (Fact \ref{Gabber}, (\ref{Summary}.3) and (\ref{Pr}.1)) and by definition $[H_{\ell^\infty}]$ maps to $[H_{\ell^\infty \Q_\ell}]$ \textit{via} $$H^1(\Pi_{\ell^\infty}, A_{\ell^\infty}\otimes B_{\ell^\infty}^{\vee})\rightarrow H^1(\Pi_{\ell^\infty}, A_{\ell^\infty \Q_\ell}\otimes B_{\ell^\infty \Q_\ell}^\vee).$$
 \noindent By (\ref{Summary}.2.2),  $[H_{\ell^\infty \Q_\ell}]=0$. So it is enough to show that $$H^1(\Pi_{\ell^\infty}, A_{\ell^\infty}\otimes B_{\ell^\infty}^{\vee})\rightarrow H^1(\Pi_{\ell^\infty}, A_{\ell^\infty \Q_\ell}\otimes B_{\ell^\infty \Q_\ell}^\vee)$$
 is injective that is $H^1(\Pi_{\ell^\infty}, A_{\ell^\infty}\otimes B_{\ell^\infty}^\vee)[\ell] = 0$. But this follows from (\ref{Summary}.3.2) applied to the $\Pi_{0\ell^\infty}$-equivariant embedding with torsion-free cokernel
  $A_{\ell^\infty}\otimes B_{\ell^\infty}^\vee\hookrightarrow H_{\ell^\infty}\otimes H_{\ell^\infty}^\vee$.   \\
 
 \section{The Grothendieck-Serre-Tate conjectures with $\F_\ell$-coefficients}\label{Conj:Tate}
One may ask whether it is reasonable to expect the (arithmetic) positive characteristic variant of the  Grothendieck-Serre-Tate conjectures with $\F_\ell$-coefficients to hold for $\ell\gg 0$. We retain the notation and conventions of  the introduction. Let $\Lambda_\ell$ denote $\Z_\ell$, $\F_\ell$ or $\Q_\ell$. Let $K_0$ be a field finitely generated over $\F_p$ and $Y_0$ a smooth proper scheme  of dimension $d$ over $K_0$. Let $Y$ denote the base change of $Y_0$ to an algebraic closure $K$ of $K_0$. For every integer $w\geq 0$, consider the following statements.

\begin{itemize}[leftmargin=* ,parsep=0cm,itemsep=0cm,topsep=0cm]
\item (\ref{Conj:Tate}.1, $\Lambda_\ell$, $w$)  (semisimplicity) The action of $\pi_1(K_0)$ on $\SH^{2w}(Y ,\Lambda_\ell)$ is semisimple. 
\item (\ref{Conj:Tate}.2, $\Lambda_\ell$, $w$)  (fullness) The map 
$Z^w(Y_0 )\otimes\Lambda_\ell\rightarrow \SH^{2w}(Y ,\Lambda_\ell(w))^{\pi_1(K_0)}$
is surjective. Here $Z^w(Y_0 )$ denotes the $\Z$-module of codimension $w$ algebraic cycles.\\
\end{itemize}

\noindent Theorem \ref{Th:GSS} and  Theorem \ref{Th:TensorInvariants} show

\begin{corollaire}\label{Final} The assertions  (\ref{Conj:Tate}.1, $\Q_\ell$, $w$), (\ref{Conj:Tate}.2, $\Q_\ell$, $i$), $i=w,d-w$ for  $Y_0$  imply  the assertions (\ref{Conj:Tate}.1, $\F_\ell$, $w$), (\ref{Conj:Tate}.2, $\F_\ell$, $w$) for  $Y_0$  provided $\ell\gg 0$ (depending on $Y_0$).
\end{corollaire}
\noindent\textit{Proof.} Assume (\ref{Conj:Tate}.1, $\Q_\ell$, $w$), (\ref{Conj:Tate}.2, $\Q_\ell$, $i$), $i=w,d-w$  for  $Y_0$. From \cite[Lem. 3.1]{MR}, for $\ell\gg 0$ (depending on $Y_0$), the canonical morphism $$Z^w(Y_0)\otimes\Z_\ell\rightarrow \SH^{2w}(Y ,\Z_\ell(w))^{\pi_1(K_0)}$$
is surjective hence, in particular, $Z^w(Y_0)\otimes\Z_\ell\rightarrow \SH^{2w}(Y ,\Z_\ell(w))$ has torsion-free cokernel. This, together with Fact \ref{Gabber}, shows that the images of $Z^w(Y_0 )\otimes\Lambda_\ell\rightarrow \SH^{2w}(Y,\Lambda_\ell(w))^{\pi_1(K_0)}$ for $\Lambda_\ell=\Z_\ell$, $\F_\ell$ or $\Q_\ell$ have the same rank. 
Thus it is enough to show (\ref{Conj:Tate}.1, $\F_\ell$, $w$) and  
$$\begin{tabular}[t]{ll} 
 (\ref{Conj:Tate}.2', $\F_\ell$, $w$) &$\SH^{2w}(Y,\Lambda_\ell(w))^{\pi_1(K_0)}$ for $\Lambda_\ell=\F_\ell,\Q_\ell$ have the same dimension for $\ell\gg 0$.\\
\end{tabular}$$
\indent From Lemma \ref{Lem:Localization}, one may freely replace $K_0$ by a finite Galois extension. In particular, we may assume that the Zariski closure of the image of $\pi_1(K_0)$ acting on $\SH^{2w}(Y,\Q_\ell(w))$ is connected for every prime $\ell\not=p$ (see Fact \ref{DeligneLP}).\\
\indent  If $K_0$ is a finite field,  (\ref{Conj:Tate}.1, $\F_\ell$, $w$) follows from the fact that the minimal   polynomial of the Frobenius acting on $\SH^w(Y ,\Q_\ell)$ is in $\Q[T]$, separable and independent of $\ell$ while 
 (\ref{Conj:Tate}.2', $\F_\ell$, $w$)  follows  from the fact that (under (\ref{Conj:Tate}.1, $\Lambda_\ell$, $w$)) the dimension of $\SH^{2w}(Y,\Lambda_\ell(w))^{\pi_1(K_0)}$ is the multiplicity of $1$  among the roots of the characteristic polynomial of Frobenius, which is in $\Q[T]$ and independent of $\ell$. \\
 \indent If $K_0$ is finitely generated, $Y_0$ is the generic fiber of a smooth proper morphism $\mathcal{Y}_0\rightarrow \mathcal{X}_0$ with $\mathcal{X}_0$ a smooth, separated, geometrically connected scheme over a finite field $k_0$. Let $\mathcal{X}$ denote the base-change of $\mathcal{X}_0$ to the algebraic closure $k$ of $k_0$. Up to enlarging $k_0$, we can find $x_0\in \mathcal{X}_0(k_0)$ such that the Frobenius $F_{x_0}$ acts semisimply on $\SH^w((\mathcal{Y}_0)_x,\Q_\ell)\simeq \SH^w(Y,\Q_\ell)$ for every $\ell\neq p$ (see the second paragraph after Prop. 1.1 in \cite{LarsenPinkAV}). Here $x$ is any geometric point over $x_0$.
 By the above argument, $F_{x_0}$ acts semisimply on $\SH^w(Y,\F_\ell)$ for $\ell\gg 0$. In particular its image is of prime-to-$\ell$ order. Thus  (\ref{Conj:Tate}.1, $\F_\ell$, $w$) follows from Theorem \ref{Th:GSS}  and \cite[Lem. 5 (b)]{SerreSS}. For (\ref{Conj:Tate}.2', $\F_\ell$, $w$), observe that
$$\SH^{2w}(Y,\Lambda_\ell(w))^{\pi_1(K_0)}=(\SH^{2w}(Y,\Lambda_\ell(w))^{\pi_1(\mathcal{X},x)})^{F_{x_0}}.$$ From Theorem \ref{Th:TensorInvariants}, $\SH^{2w}(Y,\Lambda_\ell(w))^{\pi_1(\mathcal{X},x)}$ for $\Lambda_\ell=\F_\ell,\Q_\ell$ have the same dimension. Thus  (\ref{Conj:Tate}.2', $\F_\ell$, $w$) follows from the fact that (under (\ref{Conj:Tate}.1, $\Lambda_\ell$, $w$)) the characteristic polynomial of $F_{x_0}$ acting on $\SH^{2w}(Y,\Lambda_\ell(w))^{\pi_1(\mathcal{X},x)}$ is in $\Q[T]$ and independent of $\ell$ (see (the proof of) \cite[Prop. 2.1]{LarsenPinkAV}). $\square$\\

\noindent \begin{tabular}[t]{l}
\textit{anna.cadoret@polytechnique.edu}\\
Centre de Math\'ematiques Laurent Schwartz -- \'Ecole Polytechnique\\
\end{tabular}\\

\noindent \begin{tabular}[t]{l}
\textit{pslnfq@gmail.com}\\
Department of Mathematics -- Vrije Universiteit Amsterdam \\
\end{tabular}\\

\noindent \begin{tabular}[t]{l}\textit{tamagawa@kurims.kyoto-u.ac.jp}\\
Research Institute for Mathematical Sciences -- Kyoto University \\
\end{tabular}\\

  \end{document}